\documentclass[3p,times]{elsarticle}

\usepackage{graphicx}
\usepackage{booktabs}
\usepackage{tablefootnote}
\usepackage{amsfonts}
\usepackage{amsmath}
\usepackage{amsthm}
\usepackage{amssymb}

\usepackage{algorithm}
\usepackage{caption} %caption居中
\usepackage{algorithmic}

\usepackage{geometry}
\usepackage[colorlinks,
linkcolor=black,
anchorcolor=blue,
citecolor=black]{hyperref} %%% Reference hyperlink
\makeatletter
\def\ps@pprintTitle{%
  \let\@oddhead\@empty \let\@evenhead\@empty
  \def\@oddfoot{\hfil\thepage\hfil}\let\@evenfoot\@oddfoot}
\makeatother
\newtheorem{theorem}{Theorem}[section]
\newtheorem{lemma}{Lemma}[section]
\newtheorem{proposition}{Proposition}[section]
\newtheorem{assumption}{Assumption}[section]
\newtheorem{remark}{Remark}[section]
\newtheorem{Definition}{Definition}[section]

\usepackage{bm}
\numberwithin{equation}{section}
\newcommand{\inner}[2]{\left\langle #1,\; #2\right\rangle}

\usepackage{float}
\journal{}

\begin{document}

\begin{frontmatter}

\title{Bregman Proximal Dogleg Opportunistic Majorization with Extrapolation and Structured Hybrid Barzilai--Borwein Steplengths for Nonconvex Sparse Quadratic Optimization}
 
\author[mymainaddress1,mymainaddress3]{Feifei Zhao}
\ead{zhaofeifei@shzu.edu.cn}
\author[mymainaddress1]{Qingsong Wang}
\ead{nothing2wang@hotmail.com}
\author[mymainaddress1]{Jibiao Yuan}
\ead{yuan2112020@163.com}
\author[mymainaddress1]{Zheng Peng\corref{mycorrespondingauthor}}
\cortext[mycorrespondingauthor]{Corresponding author}
\ead{pzheng@xtu.edu.cn}
\address[mymainaddress1]{School of Mathematics and Computational Science, Xiangtan University, Xiangtan 411105, P.R. China}
\address[mymainaddress3]{College of Sciences, Shihezi University, Xiangyang Street, Shihezi, 832003, Xinjiang, China}

%% Abstract
\begin{abstract}
This paper considers the minimization of the sum of a smooth quadratic function and a nonsmooth, nonconvex term whose
proximal mapping is tractable. Building on proximal dogleg opportunistic majorization with extrapolation (PDOME), we propose SHBB-PDOME, which combines a structured
Barzilai--Borwein (BB) steplength strategy with a Bregman proximal mapping. The fixed steplength used to construct the dogleg path is replaced by an adaptive steplength obtained
by blending a structured BB steplength with the local steplength stored at the preceding iteration and clipping the result. A fixed positive definite diagonal metric defines
both the weighted BB formulas and the quadratic Bregman kernel,
allowing coordinate-dependent proximal scaling. Under suitable assumptions, we establish sufficient decrease
of a Lyapunov function and show that every cluster point in the interior of the kernel domain is critical. Under additional assumptions, including the
Kurdyka--{\L}ojasiewicz (K{\L}) property, the entire sequence has finite length and converges to a critical point.
An additional relative error condition yields convergence rates determined by the K{\L} exponent. Numerical experiments on nonconvex sparse recovery, sparse approximation, and robust principal component analysis
illustrate the effectiveness of the proposed method relative to the tested alternatives.
\end{abstract}

%% Keywords
\begin{keyword}
Nonconvex sparse optimization; Structured Barzilai--Borwein steps; Opportunistic majorization; Bregman proximal methods.
\end{keyword}

\end{frontmatter}

\section{Introduction}
\label{sec1}
%% Labels are used to cross-reference an item using \ref command.
We consider composite optimization problems in which a convex quadratic loss is combined with a nonsmooth, nonconvex term:
\begin{equation}
\label{eqz2}
\min_{{x}\in\mathbb{R}^n}Q({x}):=s({x})+r({x}),
\end{equation}
where ${x}$ is the decision variable, $s$ is a quadratic function with $\nabla^2 s({x})\succcurlyeq 0$\tablefootnote{Throughout, $\mathbf M:=\nabla^2 s\succeq 0$ denotes the Hessian of the quadratic term $s$ and may be singular. To construct the Newton endpoint of the dogleg path, we add $\beta I$ to $\mathbf M$, with $\beta>0$, as described in Subsection~2.2. This gives the positive definite regularized Hessian $\mathbf B_\beta=\mathbf M+\beta I\succ 0$. This regularization is used to construct the Newton endpoint and leaves the original objective $Q$ unchanged.} and
$ r: \mathbb{R}^n \to (-\infty,+\infty]$ is a proper, lower semicontinuous, nonconvex and nonsmooth function (which may represent regularization, constraints, or complex structures). We further assume that $Q$ is lower bounded, i.e., there exists a real number $\underline{Q}$ such that $\forall {x} \in \mathbb{R}^n$, $Q({x})\geq \underline{Q}$. Additionally, $r$ is proximal-friendly, which means that the proximal operator \[
\text{prox}_{\tau r}(y) := \arg\min_{x \in \mathbb{R}^n} \left\{ r(x) + \frac{1}{2\tau} \| x - y \|^2 \right\}
\] (with step size $\tau > 0$) is easy to compute \cite{combettes2011proximal}.

Problem~\eqref{eqz2} covers a range of applications in signal processing and machine learning, including sparse signal recovery \cite{ghayem2017sparse}, audio processing \cite{stankovic2018analysis}, channel estimation \cite{sharp2008application}, feature extraction \cite{xu2025improved}, and blind source separation \cite{zibulevsky2001blind}. In sparse models, the regularizer may be convex, such as the $\ell_1$ norm \cite{candes2008introduction}, or nonconvex, such as the $\ell_{1/2}$ quasi-norm \cite{xu2012l_}, the $\ell_0$ pseudo-norm \cite{blumensath2009iterative}, and the capped $\ell_1$ penalty \cite{zhang2010analysis}. Related composite formulations also arise in low-rank matrix completion \cite{huang2021robust} and robust principal component analysis \cite{candes2011robust}.

\begin{table}[htbp]
\centering
\caption{Comparison of the methods reviewed in this section in terms of second-order information (SOI), opportunistic majorization (OM), a structured BB step size, the distance used in the proximal subproblem, and the assumption made on $r$ together with the convergence guarantee\tablefootnote{``Global'' indicates that criticality of cluster points and a convergence result for the generated sequence are available; ``Local'' indicates that only a local analysis is given.}. }
\label{tab:algorithms}
\renewcommand{\arraystretch}{1.1}

\resizebox{\textwidth}{!}{
\begin{tabular}{ccccccc}
\toprule
\textbf{Algorithm} & \textbf{Convexity of $r$} & \textbf{SOI} & \textbf{OM Strategy} & \textbf{Structured BB step} & \textbf{Bregman distance} & \textbf{Convergence} \\
\midrule
PG \cite{combettes2011proximal}                 & Yes & No & No & No & No & Global \\
PN \cite{lee2014proximal}                  & Yes & Yes & No & No & No & Local \\
mAPG \cite{li2015accelerated}                &Yes  & No & No & No & No & Global  \\
APGnc \cite{li2017convergence}                & Yes & No & No & No & No & Local \\
PANOC \cite{stella2017simple}            & No & Yes  & No & No & No & Global \\
PANOC+ \cite{de2022proximal}           & No &Yes  & No & No & No & Global \\
PDOM \cite{zhou2024proximal}             & No & Yes & Yes & No & No & Local \\
PDOME \cite{zhao2025accelerated}    & No & Yes & Yes & No & No & Global \\
\textbf{SHBB-PDOME (ours)}     & No & Yes & Yes & \textbf{Yes} & \textbf{Yes} & Global \\
\bottomrule
\end{tabular}
}
\end{table}

Various algorithms exist for solving \eqref{eqz2}, among which the proximal gradient (PG) method \cite{combettes2011proximal} is a widely used method due to its simplicity and ease of implementation. The PG method iteratively combines gradient descent on the smooth component with a proximal update for the nonsmooth term. However, in nonconvex settings, the PG method suffers from slow convergence \cite{attouch2009convergence, attouch2013convergence}.
To accelerate convergence, inertial and extrapolation techniques (e.g., the accelerated proximal gradient method for nonconvex programming, APGnc \cite{li2017convergence}, and the accelerated proximal gradient (APG) method \cite{li2015accelerated, beck2009fast}) have been developed. These methods employ adaptive mechanisms to dynamically select between standard PG updates and accelerated variants based on objective function values.
Although these accelerated variants offer improvements, their guarantees of faster convergence often require specific problem structures or assumptions \cite{li2015accelerated, attouch2009convergence, attouch2013convergence}. Consequently, the desired acceleration may fail to materialize in more general or challenging nonconvex optimization scenarios. Adaptive stepsize strategies provide another approach
to improving PG methods.
For example, Ye et al.~\cite{ye2025adaptive} developed a line-search-free proximal gradient method based on local curvature estimates, although its convergence
guarantees assume a convex nonsmooth term and thus do not directly cover the nonconvex regularizers considered here.

While adaptive PG methods use local curvature estimates to select scalar stepsizes, Newton-type methods incorporate curvature information into the proximal
subproblem itself. A representative example is the proximal Newton (PN) method \cite{lee2014proximal}, which has attracted
considerable attention for minimizing composite objectives comprising a twice differentiable term and a proximal-friendly term. Each iteration solves a proximal subproblem defined
by a Hessian-based scaled proximal operator (SPO). For convex objectives, PN achieves local superlinear convergence under suitable regularity and subproblem accuracy conditions. However, its extension to nonconvex problems faces two challenges. First, proximal subproblems can be costly to solve. Although efficient proximal trust-region subproblem solvers have been developed for smooth nonconvex objectives with nonsmooth convex terms \cite{baraldi2025efficient}, they do not provide a general, inexpensive procedure for evaluating the Hessian-based SPO. Second, the above convergence guarantees do not directly extend to nonconvex settings. These challenges motivate methods that exploit second-order information with affordable subproblem computations and convergence guarantees for nonconvex problems.  

Quasi-Newton approaches
\cite{stella2017simple, de2022proximal, patrinos2013proximal,
stella2017forward, themelis2018forward}
avoid directly solving the scaled proximal subproblem
by instead working with the forward-backward envelope (FBE).
Under suitable conditions, the FBE shares the same local
minimizers as the original objective.
These methods construct search directions using BFGS or
L-BFGS updates \cite{liu1989limited} and perform a line search based on the FBE.
Maintaining the quasi-Newton approximation requires
storing and updating a matrix or a collection of history
vectors, with computational and memory costs that grow
with the problem dimension.
Moreover, the approximation need not exploit the
structure of the exact Hessian.
The Proximal Averaged Newton-type method for Optimal
Control (PANOC) \cite{stella2017simple} combines forward-backward steps with an FBE-based line search and achieves superlinear convergence under suitable
assumptions. PANOC+ \cite{de2022proximal} introduces an adaptive
step-size rule tailored to the proximal-gradient oracle. Its analysis also addresses robustness with respect to inexact subproblem solutions. Beyond these FBE-based approaches, Ouyang and Milzarek \cite{ouyang2025trust} reformulate the stationarity conditions as a normal map-based root-finding problem and embed inexact semismooth Newton steps within a trust region framework. Under suitable assumptions, a new merit function and acceptance mechanism ensure global convergence and a transition to local (q)-superlinear convergence. The framework is also compatible with the KŁ inequality, allowing a refined convergence analysis. 

A different approach is taken by the proximal dogleg opportunistic majorization (PDOM) algorithm \cite{zhou2024proximal}, which adopts a majorization-minimization (MM) framework for problem \eqref{eqz2}.
It constructs a surrogate model along a dogleg path
formed from gradient and Newton directions and combines
this construction with a proximal mapping.
Since the smooth term is quadratic, its Hessian remains
constant across iterations.
The inverse of the Hessian, or of a fixed regularized
Hessian when needed, can therefore be precomputed
and reused.
The convergence analysis of PDOM includes a local
convergence rate result.
Subsequent work \cite{zhao2025accelerated} introduced
PDOME by incorporating extrapolation and an improved
backtracking line search into the PDOM framework.
Under suitable assumptions, every cluster point of
the generated sequence is shown to be critical.
PDOME uses a fixed steplength in its gradient branch
and a Euclidean distance in its proximal subproblem.
Building on PDOME, we incorporate a structured hybrid
BB steplength strategy and a Bregman proximal mapping
induced by a diagonal metric constructed from the data.
The metric is used both to define the weighted inner
products in the structured BB formula and to provide
coordinatewise scaling in the proximal subproblem.
The method retains the PDOME framework of extrapolation,
dogleg search, and comparison with a safeguarded
proximal candidate.
Table~\ref{tab:algorithms} summarizes the main features
of the methods discussed above.
The main contributions are as follows.

\begin{itemize}
\item We develop a structured hybrid BB steplength
strategy within a Bregman proximal dogleg framework.
The structured BB formula uses inner products weighted
by a data-dependent diagonal metric.
The resulting BB steplength is directly blended
with the local steplength stored from the previous
iteration.
The same metric induces the Bregman distance
in the proximal subproblem, providing coordinatewise
scaling.
For the separable regularizers considered in the
experiments, the resulting proximal subproblems
can be solved coordinatewise.

\item Under suitable assumptions, we establish
sufficient decrease of a Lyapunov function and prove
that every cluster point of the generated sequence
is a critical point.
If, in addition, the Lyapunov function satisfies the
K{\L} property, the iterates are bounded, the sum of
the norms of the gradient mismatches is finite, and
the ratios of consecutive step norms are uniformly
bounded, then the entire sequence has finite length
and converges to a critical point.
The convergence rate result further requires that
the norm of the gradient mismatch be controlled by
the norm of the corresponding step.

\item We evaluate SHBB-PDOME on nonconvex sparse recovery
with $\ell_0$ regularization, sparse approximation with
$\ell_{1/2}$ regularization, and robust principal component
analysis. We compare it with representative optimization methods
in terms of recovery accuracy, recovery success rates, and iteration counts. 
\end{itemize}

The rest of this paper is structured as follows: Section \ref{sec:preliminaries} presents mathematical preliminaries and related preparatory work; Section \ref{sec:3} develops the proposed SHBB-PDOME algorithm and its convergence analysis; Section \ref{sec:numerical_experiment} reports relevant experimental results; Section \ref{sec:5}  summarizes the paper and outlines future research directions.

\section{Preliminaries}\label{sec:preliminaries}
Throughout this paper, $\mathbb{R}^n$ denotes
the $n$-dimensional Euclidean space,
$\langle\cdot,\cdot\rangle$ denotes its standard
inner product, and the superscript $\top$
denotes transpose.
For $x\in\mathbb{R}^n$, define
$\|x\|:=\sqrt{x^\top x}$,
$\|x\|_1:=\sum_{i=1}^n|x_i|$, and
$\|x\|_0:=|\operatorname{supp}(x)|$,
where $\operatorname{supp}(x):=\{i:x_i\ne0\}$
is the support of $x$ and $|\cdot|$ denotes
the cardinality of a set.
Here, $\|\cdot\|$ is the Euclidean norm,
$\|\cdot\|_1$ is the $\ell_1$ norm, and
$\|\cdot\|_0$ is the $\ell_0$ pseudo-norm.
For a symmetric positive semidefinite matrix
$\mathbf M\in\mathbb{R}^{n\times n}$,
define the weighted seminorm
$\|x\|_{\mathbf M}:=\sqrt{x^\top\mathbf M x}$;
it is a norm when $\mathbf M$ is positive definite.
Given a closed set $\Omega \subseteq \mathbb{R}^n$, $\operatorname{dist}(x, \Omega) := \inf \{\|y - x\|_2 : y \in \Omega\}$ calculates the distance between $x$ and $\Omega$.

\begin{Definition}[Properness and lower semicontinuity
\cite{bertsekas2003convex}]\label{def2.1}
A function $f:\mathbb{R}^n\to(-\infty,+\infty]$
is called proper if
$\operatorname{dom}f
:=\{x\in\mathbb{R}^n:f(x)<+\infty\}$
is nonempty.
It is lower semicontinuous at $x_0\in\mathbb{R}^n$
if
\[
f(x_0)\le\liminf_{x\to x_0}f(x).
\]
It is lower semicontinuous if this condition
holds at every point of $\mathbb{R}^n$.
\end{Definition}
%\end{sloppypar}

\begin{Definition} \label{def2.2} (Gradient Lipschitz continuity \cite{bertsekas2003convex}). Let $L_f\geq 0$. A differentiable function $f:\mathbb{R}^n\to\mathbb{R}$ is said to have a Lipschitz-continuous gradient if for all 
$x,y\in\operatorname{dom}f$ it holds that
\begin{align}
\|\nabla f({x})-\nabla f({y})\|\leq L_f\|{x}-{y}\|.\nonumber
\end{align}
The corresponding Lipschitz constant is denoted as $\displaystyle L_f$.
\end{Definition}

\begin{Definition} \label{def2.3} (Subdifferential     \cite{tyrrell1998variational}). Let 
$ f:\mathbb{R}^{n}\to\mathbb{R}\cup\{+\infty\}$ be a proper and lower semicontinuous function. For a given 
$ x\in\operatorname{dom}f$, the $Fr\acute{e}chet$ subdifferential of $f$ at $x$, written as
$\hat{\partial}f({x})$, is the set of all vectors ${v}\in\mathbb{R}^n$ which satisfy
$$\liminf_{y \to x, y \neq x} \frac{f(y) - f(x) - \langle v, y - x \rangle}{\|y - x\|} \geq 0.$$
For $x\notin\operatorname{dom}f$, we set
$\hat{\partial}f({x}):=\emptyset$.
The subdifferential (which is also called the limiting subdifferential) of $f$ at 
$x\in\mathbb{R}^{n}$, written as $\partial f(x)$, is defined by
\begin{align}
\partial f({x}):=\{{v}\in\mathbb{R}^n:\exists{x}^k\to{x}, f\left({x}^k\right)\to f({x}),{v}^k\in\hat{\partial}f\left({x}^k\right)\to{v}, k\to\infty\}.\nonumber
\end{align}
As before $\partial f(x):=\emptyset$ for
$x\notin\operatorname{dom}f$, and its domain is 
$\operatorname{dom}~\partial f:=\{x\in \mathbb{R}^n:\partial f(x)\neq\emptyset\}$.
\end{Definition}

\begin{Definition}\label{def2.4} (K{\L} property \cite{attouch2010proximal}). Let $f:\mathbb{R}^n\to\mathbb{R}\cup\{+\infty\}$ be proper and closed. Then $f$ is said to have the K{\L} property at $\hat{x}\in\operatorname{dom} \partial f$ if there exist $\varsigma\in(0,+\infty]$, a neighborhood $\mathcal{B}_{\varrho}(\hat{x})\triangleq\{x:\|x-\hat{x}\|<\varrho\}$, and a continuous desingularizing concave function $\psi:[0,\varsigma)\to[0,+\infty)$ with $\psi(0)=0$ such that\\
(i) $\psi$ is continuously differentiable with $\psi^{\prime}(x)>0$ for all $x\in(0,\varsigma)$,\\
(ii) for every $x\in\mathcal B_{\varrho}(\hat{x})$
satisfying
$f(\hat{x})<f(x)<f(\hat{x})+\varsigma$, it holds that
\begin{align}
\psi'(f({x})-f(\hat{{x}}))\operatorname{dist}(0,\partial f({x}))\geq1.\nonumber
\end{align}
A proper closed function $f$ satisfying the K{\L} property at all points in $\operatorname{dom}\partial f$ is called a K{\L} function; the class of corresponding desingularizing functions is denoted by $Y_{\varsigma}$.
\end{Definition}

Every proper closed semi-algebraic function
is a K{\L} function \cite{attouch2010proximal}.
In particular, the composite objectives considered
in Section~\ref{sec:numerical_experiment}
are semi-algebraic.
The finite-length and convergence results below
do not require a prescribed K{\L} exponent,
whereas the convergence-rate analysis uses
a power-type desingularizing function.

\begin{lemma}[Uniformized K{\L} property
\cite{bolte2014proximal}]\label{lem:2.6}
Let $\mathbb W\subseteq\mathbb R^n$ be a nonempty
compact set, and let
$f:\mathbb R^n\to\mathbb R\cup\{+\infty\}$
be proper and lower semicontinuous.
Assume that $f$ is constant on $\mathbb W$
and satisfies the K{\L} property at every point
of $\mathbb W$.
Then there exist $\varepsilon>0$,
$\varsigma>0$, and $\psi\in Y_{\varsigma}$
such that
\[
\psi'\bigl(f(z)-f(\bar z)\bigr)
\operatorname{dist}(0,\partial f(z))\ge1
\]
for every $\bar z\in\mathbb W$ and every
$z\in\mathbb R^n$ satisfying
$\operatorname{dist}(z,\mathbb W)<\varepsilon$
and
$f(\bar z)<f(z)<f(\bar z)+\varsigma$.
\end{lemma}

\subsection{The Proximal Gradient Method from an MM Perspective}
\label{sec:pg_mm}
The PG method with extrapolation admits a
majorization-based interpretation.
At each iteration, it minimizes, with respect to $x$,
the surrogate
\begin{equation}\label{mpg}
m_{\mathrm{pg}}(x;y^k)+r(x),\qquad
m_{\mathrm{pg}}(x;y^k)
:=s(y^k)+\inner{g_k}{x-y^k}
+\frac{1}{2\eta}\|x-y^k\|^2,
\end{equation}
where $y^k=x^k+\zeta(x^k-x^{k-1})$
is the extrapolation point and
$g_k=\nabla s(y^k)$.
The subproblem can be solved using the proximal
mapping of $r$.
Extrapolation and related inertial techniques
incorporate information from previous iterates
and are used in several proximal methods
\cite{ochs2014ipiano,yang2024proximal,wen2017linear}.

For $0<\eta<1/L_s$, the model satisfies
$m_{\mathrm{pg}}(x;y^k)\ge s(x)$
for all $x\in\mathbb R^n$.
The fixed step size $\eta$ used here is also the
safeguard step of Algorithm~\ref{Algorithm3}; it is
distinct from the local steplength $\eta_\mu$, which
is determined adaptively along the dogleg direction
for the current $\mu$.
The MM interpretation of the PG method and the
corresponding majorization relation are discussed
in \cite[Section~2.1]{zhao2025accelerated}.
In the presence of extrapolation, this majorization
property alone does not guarantee monotonic decrease
of the objective values.
Sufficient decrease of a Lyapunov function
for the proposed method is established under
suitable parameter conditions in
Lemma~\ref{lem:lyap}.

The gradient step uses only first-order information
about the smooth term, and gradient methods can
converge slowly on ill-conditioned problems
\cite[Chapter~9]{boyd2004convex}.
This motivates the dogleg construction introduced
next, which combines gradient and Newton directions
within the majorization framework.
\subsection{Bregman Geometry and the Hybrid Dogleg Direction}
\label{sec:hybrid}
We build on the PDOME algorithm
\cite{zhao2025accelerated}, which incorporates
a dogleg search inspired by trust-region methods
into the opportunistic majorization framework.
We introduce an adaptive structured hybrid BB
steplength strategy to determine the scale
of the gradient branch of the dogleg path.
We also replace the Euclidean quadratic term
in the proximal surrogate with a suitably scaled
Bregman divergence $D_{\hat h}$.
We first introduce the general Bregman geometry
and then specialize it to a quadratic kernel
induced by a data-dependent diagonal metric.
The same metric defines the weighted inner products
in the structured BB formulas, linking
the steplength construction with the geometry
of the proximal subproblem.

\noindent\textbf{Bregman geometry.}
Let $\hat h:\mathbb R^n\to(-\infty,+\infty]$ be
proper, closed, convex, and essentially smooth, with
$\operatorname{int}\operatorname{dom}\hat h
\neq\varnothing$.
The Bregman divergence generated by $\hat h$
\cite{bregman1967relaxation} is
\begin{equation}\label{eq:3pt}
D_{\hat h}(x,y)
:=\hat h(x)-\hat h(y)
-\langle\nabla\hat h(y),x-y\rangle,
\qquad
x\in\operatorname{dom}\hat h,\quad
y\in\operatorname{int}\operatorname{dom}\hat h.
\end{equation}
If $\hat h$ is twice continuously differentiable
with a positive definite Hessian throughout
$\operatorname{int}\operatorname{dom}\hat h$,
its Hessian defines the local inner product
and associated norm:
\[
\langle u,v\rangle_{\hat h,y}
:=u^\top\nabla^2\hat h(y)v,
\qquad
\|u\|_{\hat h,y}^{2}
:=\langle u,u\rangle_{\hat h,y}.
\]
For $\hat h_0(x):=\frac12\|x\|^2$,
these coincide with the Euclidean inner product
and norm, and
$D_{\hat h_0}(x,y)=\frac12\|x-y\|^2$. We now specialize to the quadratic kernel
$h_{\mathbf M_{\mathrm d}}(x)
:=\frac12x^\top\mathbf M_{\mathrm d}x$,
where
$\mathbf M_{\mathrm d}
=\operatorname{Diag}(m_1,\ldots,m_n)$
is fixed throughout the iterations and $m_i>0$.
The associated inner product and norm satisfy
\begin{equation}\label{eq:Md-metric}
\langle u,v\rangle_{\mathbf M_{\mathrm d}}
:=u^\top\mathbf M_{\mathrm d}v,
\qquad
\|u\|_{\mathbf M_{\mathrm d}}^2
:=u^\top\mathbf M_{\mathrm d}u.
\end{equation}
Its Bregman divergence is
\begin{equation}\label{eq:Md-div}
D_{h_{\mathbf M_{\mathrm d}}}(u,v)
=\frac12\|u-v\|_{\mathbf M_{\mathrm d}}^2.
\end{equation}
For $\hat h=h_{\mathbf M_{\mathrm d}}$, we have
$\nabla^2\hat h=\mathbf M_{\mathrm d}$. Its smallest and
largest eigenvalues are
\[
\lambda_{\min}^{\hat h}:=\min_{1\le i\le n}m_i>0,
\qquad
\lambda_{\max}^{\hat h}:=\max_{1\le i\le n}m_i.
\]
We use the same metric $\mathbf M_{\mathrm d}$ to construct the
weighted BB steplengths below.

\noindent\textbf{Structured BB steplengths.}
For $k\ge1$, let $\Delta_k:=x^k-x^{k-1}$.
Write the least-squares loss as
$s(x)=\frac12\|R(x)\|^2$, where
$R:\mathbb R^n\to\mathbb R^m$ is the residual mapping
and $J(x)\in\mathbb R^{m\times n}$ is its Jacobian.
Set $R_k:=R(x^k)$ and $J_k:=J(x^k)$.
For $\Delta_k\neq0$, we follow the structured
curvature construction in \cite{awwal2022barzilai}
and define
\[
\vartheta_{k-1}
:=3R_k^\top
\big[(J_k-J_{k-1})\Delta_k
-2(R_k-R_{k-1})\big],
\]
\[
\bar y_{k-1}
:=J_k^\top J_k\Delta_k
+(J_k-J_{k-1})^\top R_k
+\frac{\vartheta_{k-1}}{\|\Delta_k\|^2}\Delta_k.
\]
Although the curvature formula is presented using general residual notation, we use it here only
for quadratic least-squares problems.

Let $\varepsilon_{\mathrm{BB}}>0$.
Whenever
$\langle\Delta_k,\bar y_{k-1}\rangle_{\mathbf M_{\mathrm d}}
>\varepsilon_{\mathrm{BB}}$,
define the weighted BB steplengths by
\begin{equation}\label{eq:bbstep}
\tau_{\mathrm{BB1}}^{\mathbf M_{\mathrm d}}
:=\frac{\|\Delta_k\|_{\mathbf M_{\mathrm d}}^2}
{\langle\Delta_k,\bar y_{k-1}\rangle_{\mathbf M_{\mathrm d}}},
\qquad
\tau_{\mathrm{BB2}}^{\mathbf M_{\mathrm d}}
:=\frac{\langle\Delta_k,\bar y_{k-1}\rangle_{\mathbf M_{\mathrm d}}}
{\|\bar y_{k-1}\|_{\mathbf M_{\mathrm d}}^2}.
\end{equation}
The curvature condition and
$\mathbf M_{\mathrm d}\succ0$ ensure that both
steplengths are well defined and positive.
Their convex combination is
\[
\tau_{\mathrm{SHBB}}^{\mathbf M_{\mathrm d}}
:=w\tau_{\mathrm{BB1}}^{\mathbf M_{\mathrm d}}
+(1-w)\tau_{\mathrm{BB2}}^{\mathbf M_{\mathrm d}},
\qquad w\in[0,1].
\]
At \(k=0\), we set
\(\tau_{\mathrm{SHBB}}^{\mathbf M_{\mathrm d}}=\eta\),
where \(\eta>0\).
For \(k\ge 1\), we compute
\(\tau_{\mathrm{SHBB}}^{\mathbf M_{\mathrm d}}\) as the
convex combination above if \(\Delta_k\ne 0\) and the
curvature condition holds.
Otherwise, we set
\(\tau_{\mathrm{SHBB}}^{\mathbf M_{\mathrm d}}=\tau_{k-1}\).

We directly combine the resulting BB steplength with the
local steplength stored from the previous iteration,
denoted by $\eta_{\mu^{k-1}}$. For $\rho_k\in[0,1]$, define
\begin{equation}\label{eq:tauk}
\begin{aligned}
\tau_k^{\mathrm{raw}}
&:= \rho_k\eta_{\mu^{k-1}}
  +(1-\rho_k)\tau_{\mathrm{SHBB}}^{\mathbf M_{\mathrm d}},\\[3pt]
\tau_k
&:= \max\left\{
  \tau_{\min},
  \min\left\{\tau_k^{\mathrm{raw}},\frac{1}{L_s+\beta}\right\}
  \right\}.
\end{aligned}
\end{equation}
Here $\eta_{\mu^{-1}}:=\eta$ and
$0<\tau_{\min}\le 1/(L_s+\beta)$.
Hence, by construction,
\[
\tau_{\min}\le\tau_k\le\frac1{L_s+\beta},
\qquad k\ge0.
\]
The lower bound is used in the analysis;
the upper truncation is the spectral cap.
The parameter $\rho_k$ controls the relative
contributions of the stored local steplength
and the BB steplength.
When $\rho_k=1$, the raw value
$\tau_k^{\mathrm{raw}}$ equals the stored
local steplength; when $\rho_k=0$, it equals
the BB steplength. The clipping in
\eqref{eq:tauk} is applied afterwards.
The local steplength and its storage and fallback
rules are specified in
Subsection~\ref{sec:algdev}.
Using the value stored from the previous iteration
makes $\tau_k$ available before the current
dogleg search.
The capped value $\bar\tau_k=\tau_k$ is used
to construct the hybrid dogleg path.

\noindent\textbf{Hybrid dogleg direction.}
Let $g_k:=\nabla s(y^k)$ and
$\mathbf M:=\nabla^2s\succeq0$.
Since $s$ is quadratic, $\mathbf M$ is constant.
It is distinct from the diagonal metric
$\mathbf M_{\mathrm d}$ used in the weighted
BB formulas and the Bregman proximal subproblem.
Fix $\beta>0$ and set
\begin{equation}\label{eq:precond}
\mathbf B_\beta:=\mathbf M+\beta I\succ0.
\end{equation}
The exact regularized Newton endpoint of the dogleg
path is
\begin{equation}\label{eq:newton-leg}
d_N:=-\mathbf B_\beta^{-1}g_k.
\end{equation}
The objective retains the Hessian $\mathbf M$;
$\mathbf B_\beta$ is used only to define $d_N$.
In particular, $\mathbf M$ need not be positive
definite, which covers the underdetermined case.

For the DCT-structured least-squares case used in the
experiments, $A=\Phi\Upsilon$ with $\Phi$ a row-subsampling
matrix and $\Upsilon$ the orthogonal DCT basis. Then
$\mathbf M=\Upsilon^\top\Phi^\top\Phi\Upsilon$, the matrix
$\Phi^\top\Phi$ is diagonal, and
\[
\mathbf B_\beta
=\Upsilon^\top\left(\Phi^\top\Phi+\beta I\right)\Upsilon,
\qquad
\mathbf B_\beta^{-1}g
=\Upsilon^\top\left(\Phi^\top\Phi+\beta I\right)^{-1}
\Upsilon g,
\]
so $\mathbf B_\beta^{-1}$ is computed exactly by DCT.

Since the hybrid steplength is capped at
$1/(L_s+\beta)$ in \eqref{eq:tauk}, we set
\begin{equation}\label{eq:taubar}
\bar\tau_k:=\tau_k.
\end{equation}
Then $\bar\tau_k>0$ and
$\bar\tau_k\le1/(L_s+\beta)$ by Lemma~\ref{lem:ssbb}.
The hybrid dogleg path is then defined by
\begin{equation}\label{eq:path}
d^{\mathrm{hy}}(\mu)
=
\begin{cases}
-\bar\tau_k g_k,
& \mu\in(0,1],\\[2pt]
(2-\mu)(-\bar\tau_k g_k)+(\mu-1)d_N,
& \mu\in(1,2].
\end{cases}
\end{equation}
For $\mu\in[1,2]$, the hybrid dogleg path
interpolates between the gradient step
$-\bar\tau_k g_k$ and the Newton or regularized
Newton step $d_N$.
The first branch is constant in $\mu$ and,
since Algorithm~\ref{Algorithm3} uses only
$\mu\in[1,2]$, contributes only the gradient
endpoint at $\mu=1$.
The descent property is established in
Lemma~\ref{lem:descent}.
The capped hybrid steplength $\bar\tau_k$
determines this endpoint, whereas the dogleg
proximal subproblem uses the scaled local
steplength $\gamma\eta_{\mu^k}$.
The safeguard subproblem uses the fixed
steplength $\eta$.
Subsection~\ref{sec:algdev} defines the projected
gradient and the local steplength and explains
how they determine the proximal surrogate.

For local calculations in the subsequent proofs,
we fix an iteration $k$ and use coordinates
centered at $y^k$.
Under the change of variables $z=y^k+x$,
we use $s(y^k+x)-s(y^k)$ as the translated
smooth function and retain the notation $s$.
Thus,
\[
s(x)=\langle g_k,x\rangle
+\frac12x^\top\mathbf M x,
\qquad
\nabla s(x)=g_k+\mathbf M x,
\]
where $x$ denotes displacement from $y^k$.
In these coordinates, the regularizer is
evaluated as $r(y^k+x)$, and the Bregman
kernel as $\hat h(y^k+x)$.
This change of coordinates is used only
to simplify local calculations at a fixed iteration.
All comparisons of objective values across iterations
and all Lyapunov estimates refer to the original
objective and iterates.

\section{Main Results}\label{sec:3}
In this section, we present the SHBB-PDOME algorithm
for problem~\eqref{eqz2} and establish its
convergence properties under suitable assumptions.
The algorithm retains the extrapolation,
a finite dogleg-parameter search with an acceptance
test, and the opportunistic majorization framework
of PDOME, while incorporating
the adaptive structured hybrid BB steplength strategy
and the Bregman geometry induced by a diagonal metric
introduced in Subsection~\ref{sec:hybrid}.

\subsection{Algorithm Description}\label{sec:algdev}
Using the dogleg path constructed in
Subsection~\ref{sec:hybrid}, we define the local
surrogate and the two proximal candidates.
The candidate selection follows the comparison
strategy of PDOME
\cite[Algorithm~2]{zhao2025accelerated}.
For $g_k\neq0$ and $\mu\in(0,2]$, define the
projected gradient and local steplength by
\begin{equation}
g_\mu
:=
\frac{\langle g_k,d^{\mathrm{hy}}(\mu)\rangle}
{\|d^{\mathrm{hy}}(\mu)\|^2}
d^{\mathrm{hy}}(\mu),
\qquad
\eta_\mu
:=
-\frac{\|d^{\mathrm{hy}}(\mu)\|^2}
{\langle g_k,d^{\mathrm{hy}}(\mu)\rangle}.
\label{eq:projected_gradient_local_step}
\end{equation}
Lemma~\ref{lem:descent} ensures that
$\langle g_k,d^{\mathrm{hy}}(\mu)\rangle<0$,
so both quantities are well defined and
$\eta_\mu>0$.
These definitions use the Euclidean inner product
and imply
$g_\mu=-\eta_\mu^{-1}d^{\mathrm{hy}}(\mu)$
and
\begin{equation}
\langle g_\mu,d^{\mathrm{hy}}(\mu)\rangle
=
\langle g_k,d^{\mathrm{hy}}(\mu)\rangle.
\label{eq:projected_slope_identity}
\end{equation}

For each trial parameter $\mu$, define the local
Bregman model by
\begin{equation}\label{eq:dogleg_local_model}
m_\mu^{\hat h}(x;y^k)
:=s(y^k)+\langle g_\mu,x-y^k\rangle
+\frac{1}{\eta_\mu}D_{\hat h}(x,y^k).
\end{equation}
For the quadratic kernel
$\hat h=h_{\mathbf M_{\mathrm d}}$,
the dogleg trial candidate is computed using
the scaled local steplength $\gamma\eta_\mu$:
\begin{equation}\label{eq:dogleg_bregman_candidate}
\begin{aligned}
\tilde x(\mu)
&\in\arg\min_x
\left\{
r(x)+\langle g_\mu,x-y^k\rangle
+\frac{1}{\gamma\eta_\mu}
D_{h_{\mathbf M_{\mathrm d}}}(x,y^k)
\right\}\\
&=\arg\min_x
\left\{
r(x)+\frac{1}{2\gamma\eta_\mu}
\left\|
x-\left(
y^k+\gamma\mathbf M_{\mathrm d}^{-1}
d^{\mathrm{hy}}(\mu)
\right)
\right\|_{\mathbf M_{\mathrm d}}^2
\right\}.
\end{aligned}
\end{equation}
The equality follows by completing the square,
using $g_\mu=-\eta_\mu^{-1}d^{\mathrm{hy}}(\mu)$,
and omitting a term independent of $x$. When $g_k\ne0$, we inspect the parameters in
$\mathcal T:=\{1\}\cup\{1+2^{-i}:i=0,\ldots,10\}$
in decreasing order.
For each $\mu$, we compute $\tilde x(\mu)$
and accept the parameter only if
\begin{equation}\label{eq:dogleg_acceptance}
m_\mu^{\hat h}(\tilde x(\mu);y^k)
\ge s(\tilde x(\mu)),
\qquad
\langle g_\mu-g_k,x^k-y^k\rangle\le0.
\end{equation}
The first condition requires majorization
at the trial candidate.
The second imposes a sign condition on the
Euclidean inner product of the gradient mismatch
and $x^k-y^k$.
We denote the first accepted parameter by $\mu^k$
and set $\tilde x^k:=\tilde x(\mu^k)$.
If $g_k=0$ or no parameter satisfies
\eqref{eq:dogleg_acceptance},
we omit the dogleg candidate.

At every iteration, we also compute a safeguard
candidate using the fixed steplength $\eta$:
\begin{equation}\label{eq:safeguard_bregman_candidate}
\begin{aligned}
v^{k+1}
&\in\arg\min_x
\left\{
r(x)+\langle g_k,x-y^k\rangle
+\frac{1}{\eta}
D_{h_{\mathbf M_{\mathrm d}}}(x,y^k)
\right\}\\
&=\arg\min_x
\left\{
r(x)+\frac{1}{2\eta}
\left\|
x-\left(
y^k-\eta\mathbf M_{\mathrm d}^{-1}g_k
\right)
\right\|_{\mathbf M_{\mathrm d}}^2
\right\}.
\end{aligned}
\end{equation}
The equality follows by completing the square
and omitting a term independent of $x$.
Since $r$ may be nonconvex, the proximal
subproblems need not have unique solutions.
If multiple global minimizers exist,
any one of them may be selected.

If an admissible dogleg candidate is available,
we compare its objective value with that of
the safeguard candidate and select the one
with the smaller value, choosing $\tilde x^k$
in case of a tie.
Otherwise, we set $x^{k+1}=v^{k+1}$.
The blending weight is updated according to
Algorithm~\ref{Algorithm3}.
For the next iteration, we store $\eta_{\mu^k}$
whenever an admissible dogleg candidate exists,
regardless of which candidate is selected.
If no admissible dogleg candidate exists,
we set the stored value $\eta_{\mu^k}:=\eta$.
\begin{algorithm}[t]
\caption{SHBB-PDOME: Structured Hybrid
Barzilai--Borwein Proximal Dogleg Opportunistic
Majorization with Extrapolation.}
\label{Algorithm3}
\footnotesize
\begin{algorithmic}[1]
\STATE \textbf{Input:}
$x^0\in\operatorname{dom}r$,
$\eta\in(0,\lambda_{\min}(\mathbf M_{\mathrm d})/L_s)$,
$\beta>0$,
$\varepsilon_{\mathrm{BB}}>0$,
$\gamma\in(0,1)$, $\zeta\in[0,1)$,
$w\in[0,1]$,
$\rho_0\in(0,1)$, $\rho_{\min}\in(0,\rho_0)$,
$\tau_{\min}\in(0,1/(L_s+\beta)]$,
$\delta_\rho>0$, $\mathrm{tol}>0$,
$\alpha_{\mathrm d}\ge0$, $\beta_{\mathrm d}>0$,
and $K$.

\STATE Construct $\mathbf M_{\mathrm d}\succ0$
as specified in Section~2.2.
Initialize $k=0$, $x^{-1}=x^0$,
$\tau_{-1}=\eta$, and $\eta_{\mu^{-1}}=\eta$.

\REPEAT
\STATE
$\Delta_k=x^k-x^{k-1}$,
$y^k=x^k+\zeta\Delta_k$,
$g_k=\nabla s(y^k)$.
Set
$\tau_{\mathrm{SHBB}}^{\mathbf M_{\mathrm d}}
=\tau_{k-1}$.

\IF{$\Delta_k\ne0$}
\STATE
$\vartheta_{k-1}
=3R_k^\top
\bigl[(J_k-J_{k-1})\Delta_k
-2(R_k-R_{k-1})\bigr]$,
\quad
$\displaystyle
\bar y_{k-1}
=J_k^\top J_k\Delta_k
+(J_k-J_{k-1})^\top R_k
+\frac{\vartheta_{k-1}}{\|\Delta_k\|^2}\Delta_k$.

\IF{$\langle\Delta_k,\bar y_{k-1}
\rangle_{\mathbf M_{\mathrm d}}
>\varepsilon_{\mathrm{BB}}$}
\STATE
$\displaystyle
\tau_{\mathrm{BB1}}^{\mathbf M_{\mathrm d}}
=\frac{\|\Delta_k\|_{\mathbf M_{\mathrm d}}^2}
{\langle\Delta_k,\bar y_{k-1}\rangle_{\mathbf M_{\mathrm d}}},
\quad
\tau_{\mathrm{BB2}}^{\mathbf M_{\mathrm d}}
=\frac{\langle\Delta_k,\bar y_{k-1}\rangle_{\mathbf M_{\mathrm d}}}
{\|\bar y_{k-1}\|_{\mathbf M_{\mathrm d}}^2}$,
\quad
$\tau_{\mathrm{SHBB}}^{\mathbf M_{\mathrm d}}
=w\tau_{\mathrm{BB1}}^{\mathbf M_{\mathrm d}}
+(1-w)\tau_{\mathrm{BB2}}^{\mathbf M_{\mathrm d}}$.
\ENDIF
\ENDIF

\STATE
$\tau_k^{\mathrm{raw}}
=\rho_k\eta_{\mu^{k-1}}
+(1-\rho_k)\tau_{\mathrm{SHBB}}^{\mathbf M_{\mathrm d}}$;
set
$\tau_k
=\max\{\tau_{\min},\min\{\tau_k^{\mathrm{raw}},1/(L_s+\beta)\}\}$.
Compute $\bar\tau_k$ by \eqref{eq:taubar}.

\STATE If $g_k\ne0$, compute the exact regularized
Newton step $d_N=-\mathbf B_\beta^{-1}g_k$
by \eqref{eq:newton-leg},
and search
$\mu\in\mathcal T:=\{1\}\cup
\{1+2^{-i}:i=0,\ldots,10\}$
in decreasing order.
For each $\mu$, compute
$d^{\mathrm{hy}}(\mu)$, $g_\mu$, and $\eta_\mu$
as defined above, and
$\tilde x(\mu)\in\arg\min_x
\left\{
r(x)+\langle g_\mu,x-y^k\rangle
+\frac{1}{\gamma\eta_\mu}
D_{h_{\mathbf M_{\mathrm d}}}(x,y^k)
\right\}$.
Accept the first parameter $\mu^k$ such that
$m_\mu^{h_{\mathbf M_{\mathrm d}}}
(\tilde x(\mu);y^k)\ge s(\tilde x(\mu))$
and
$\langle g_\mu-g_k,x^k-y^k\rangle\le0$,
and set $\tilde x^k=\tilde x(\mu^k)$.
If $g_k=0$ or no parameter is admissible,
skip the dogleg candidate.

\STATE
$v^{k+1}\in\arg\min_x
\left\{
r(x)+\langle g_k,x-y^k\rangle
+\eta^{-1}D_{h_{\mathbf M_{\mathrm d}}}(x,y^k)
\right\}$.

\STATE If no admissible dogleg candidate is
available, set
$x^{k+1}=v^{k+1}$,
$\rho_{k+1}
=\max\{\rho_{\min},\rho_k-\delta_\rho\}$,
and $\eta_{\mu^k}:=\eta$;
otherwise, set
\[
x^{k+1}
=
\begin{cases}
v^{k+1}, & Q(\tilde x^k)>Q(v^{k+1}),\\
\tilde x^k, & Q(\tilde x^k)\le Q(v^{k+1}),
\end{cases}
\qquad
\rho_{k+1}
=
\begin{cases}
\max\{\rho_{\min},\rho_k-\delta_\rho\},
& Q(\tilde x^k)>Q(v^{k+1}),\\
\min\{1,\rho_k+\delta_\rho\},
& Q(\tilde x^k)<Q(v^{k+1}),\\
\rho_k,
& Q(\tilde x^k)=Q(v^{k+1}).
\end{cases}
\]

\STATE Retain $\eta_{\mu^k}$ for the next iteration
and set $k\leftarrow k+1$.
\UNTIL{
$\|x^k-x^{k-1}\|/(1+\|x^k\|)<\mathrm{tol}$
or $k\ge K$
}
\STATE \textbf{Output:} $x^k$.
\end{algorithmic}
\end{algorithm}

\subsection{Convergence Analysis}
We next establish the convergence properties of SHBB-PDOME under the following assumptions.

\begin{assumption}[Problem structure and steplength parameters]
\label{ass:main}
In addition to the standing assumptions on
problem~\eqref{eqz2}, assume that the quadratic
smooth term $s$ has a positive semidefinite
Hessian $\mathbf M:=\nabla^2s\succeq0$, and let
$L_s:=\lambda_{\max}(\mathbf M)>0$.
Fix $\beta>0$ and define
$\mathbf B_\beta:=\mathbf M+\beta I\succ0$.
The Newton endpoint used in Algorithm~\ref{Algorithm3}
is $d_N=-\mathbf B_\beta^{-1}g_k$.
Since $\mathbf M\succeq0$,
$\lambda_{\min}(\mathbf B_\beta)
=\lambda_{\min}(\mathbf M)+\beta\ge\beta$ and
$\lambda_{\max}(\mathbf B_\beta)=L_s+\beta$.
The algorithmic parameters satisfy
\[
0<\tau_{\min}\le\frac1{L_s+\beta}.
\]
\end{assumption}

\begin{assumption}[Bregman geometry]\label{ass:breg}
Let $\hat h$ satisfy the kernel conditions stated
in Subsection~\ref{sec:hybrid}.
Assume further that
$\hat h\in C^2(\operatorname{int}\operatorname{dom}\hat h)$
and that there exist constants
$0<\lambda_{\min}^{\hat h}
\le\lambda_{\max}^{\hat h}<+\infty$
such that
\begin{equation}\label{eq:hmetric}
\lambda_{\min}^{\hat h}I
\preceq\nabla^2\hat h(y)
\preceq\lambda_{\max}^{\hat h}I,
\qquad
y\in\operatorname{int}\operatorname{dom}\hat h.
\end{equation}
For the diagonal quadratic kernel
$h_{\mathbf M_{\mathrm d}}$,
we take
$\lambda_{\min}^{\hat h}:=\min_i m_i$
and
$\lambda_{\max}^{\hat h}:=\max_i m_i$.
\end{assumption}

Throughout the analysis, assume that
$x^0\in\operatorname{dom}r
\cap\operatorname{int}\operatorname{dom}\hat h$
and that every proximal subproblem used by
Algorithm~\ref{Algorithm3} admits a global minimizer.
All iterates, extrapolated points, and proximal
candidates are assumed to lie in
$\operatorname{int}\operatorname{dom}\hat h$.
This domain requirement is automatic for
the quadratic kernel $h_{\mathbf M_{\mathrm d}}$.
We study the sequence generated by
Algorithm~\ref{Algorithm3}
with both stopping criteria disabled.

\begin{lemma}[Bounds on the hybrid steplength]
\label{lem:ssbb}
Suppose that Assumption~\ref{ass:main} holds.
Then, for every $k\ge0$,
\begin{equation}\label{eq:tau-bounds}
0<\tau_{\min}\le\tau_k\le\frac1{L_s+\beta}.
\end{equation}
Since $\bar\tau_k=\tau_k$ by \eqref{eq:taubar},
the same bounds hold for $\bar\tau_k$:
\begin{equation}\label{eq:hyb}
0<\tau_{\min}\le\bar\tau_k\le\frac1{L_s+\beta},
\qquad k\ge0.
\end{equation}
\end{lemma}

\begin{proof}
By \eqref{eq:tauk},
\[
\tau_k
=\max\left\{
\tau_{\min},
\min\left\{
\tau_k^{\mathrm{raw}},
\frac1{L_s+\beta}
\right\}\right\}.
\]
The lower bound is immediate. For the upper bound,
$\min\{\tau_k^{\mathrm{raw}},1/(L_s+\beta)\}\le1/(L_s+\beta)$ and
$\tau_{\min}\le1/(L_s+\beta)$ by Assumption~\ref{ass:main},
so the maximum of these two quantities is at most
$1/(L_s+\beta)$. The claim for $\bar\tau_k$ follows from
$\bar\tau_k=\tau_k$.
\end{proof}
\begin{lemma}[Descent property of the hybrid dogleg direction]
\label{lem:descent}
Suppose that Assumption~\ref{ass:main} holds
and $g_k\ne0$.
Then, for every $\mu\in(0,2]$,
\[
\langle g_k,d^{\mathrm{hy}}(\mu)\rangle<0,
\qquad
\left\langle
d^{\mathrm{hy}}(\mu),
\nabla s\bigl(d^{\mathrm{hy}}(\mu)\bigr)
\right\rangle\le0,
\]
where $s$ is expressed in the translated coordinates
centered at $y^k$.
\end{lemma}

\begin{proof}
In the translated coordinates centered at $y^k$, we have
$\nabla s(d)=g_k+\mathbf M d$. By Lemma~\ref{lem:ssbb},
$0<\bar\tau_k\le1/(L_s+\beta)\le1/L_s$.
Set $p:=-\bar\tau_k g_k$. For $\mu\in(0,1]$, the direction
$d^{\mathrm{hy}}(\mu)$ equals $p$. For $\mu\in[1,2]$, it is a
convex combination of $p$ and $d_N$.

Since $g_k\ne0$ and $\mathbf B_\beta\succ0$, we have
\[
\langle g_k,p\rangle=-\bar\tau_k\|g_k\|^2<0,
\qquad
\langle g_k,d_N\rangle=-g_k^\top\mathbf B_\beta^{-1}g_k<0.
\]
It follows that $\langle g_k,d^{\mathrm{hy}}(\mu)\rangle<0$
for every $\mu\in(0,2]$.

To prove the second inequality, define
$F(d):=\langle d,g_k+\mathbf M d\rangle$.
The function $F$ is convex because $\mathbf M\succeq0$.
At the gradient endpoint,
\[
F(p)
=-\bar\tau_k\|g_k\|^2+\bar\tau_k^2 g_k^\top\mathbf M g_k
\le \bar\tau_k(\bar\tau_k L_s-1)\|g_k\|^2
\le 0.
\]
At the regularized Newton endpoint,
$\mathbf B_\beta d_N=-g_k$ gives
$g_k+\mathbf M d_N=-\beta d_N$. Hence
\[
F(d_N)=-\beta\|d_N\|^2\le0.
\]
Therefore, $F(d^{\mathrm{hy}}(\mu))\le0$ on $(0,1]$ because
$d^{\mathrm{hy}}(\mu)=p$, and on $[1,2]$ by convexity of $F$.
This completes the proof.
\end{proof}

\begin{lemma}[Monotonicity of the local steplength in $\mu$]
\label{lem:mono}
Suppose that Assumption~\ref{ass:main} holds
and $g_k\ne0$.
Then the local steplength $\eta_\mu$ defined in
\eqref{eq:projected_gradient_local_step}
is nondecreasing in $\mu\in(0,2]$.
\end{lemma}

\begin{proof}
For $\mu\in(0,1]$, we have
$d^{\mathrm{hy}}(\mu)=-\bar\tau_k g_k$,
and hence $\eta_\mu=\bar\tau_k$.

For $\mu\in[1,2]$, let
$u:=\|g_k\|^2$,
$v:=\langle g_k,\mathbf B_\beta^{-1}g_k\rangle$,
and $q:=\|\mathbf B_\beta^{-1}g_k\|^2$.
Since $g_k\ne0$, we have $u>0$.
The spectral bound and Lemma~\ref{lem:ssbb}
give
$v\ge u/\lambda_{\max}(\mathbf B_\beta)
\ge\bar\tau_k u$,
while the Cauchy--Schwarz inequality gives
$q\ge v^2/u$.
Substituting the path definition
\eqref{eq:path} into
\eqref{eq:projected_gradient_local_step}
and rearranging, we obtain
\begin{align*}
\eta_\mu
&=
\frac{
\bar\tau_k^2(2-\mu)^2u
+2\bar\tau_k(\mu-1)(2-\mu)v
+(\mu-1)^2q
}{
\bar\tau_k(2-\mu)u+(\mu-1)v
}\\
&=
\bar\tau_k(2-\mu)+(\mu-1)\frac vu
+
\frac{
(\mu-1)^2(q-v^2/u)
}{
\bar\tau_k(2-\mu)u+(\mu-1)v
}.
\end{align*}
The denominator is strictly positive on $[1,2]$.
Differentiating with respect to $\mu\in(1,2)$
yields
\[
\frac{d\eta_\mu}{d\mu}
=
\frac vu-\bar\tau_k
+
\left(q-\frac{v^2}{u}\right)
\frac{
(\mu-1)\bigl(
2\bar\tau_k u+(\mu-1)(v-\bar\tau_k u)
\bigr)
}{
\bigl(\bar\tau_k(2-\mu)u+(\mu-1)v\bigr)^2
}
\ge0.
\]
By continuity, $\eta_\mu$ is nondecreasing
on $[1,2]$.
Since $\eta_1=\bar\tau_k$ agrees with its
constant value on $(0,1]$,
the conclusion follows.
\end{proof}

\begin{lemma}[Bounds on the local steplength]
\label{lem:eta}
Suppose that Assumption~\ref{ass:main} holds.
For every iteration $k$ with $g_k\ne0$,
the local steplength associated with the
hybrid dogleg path \eqref{eq:path} satisfies
\[
0<\eta_{\min}\le\eta_\mu\le\eta_{\max}<+\infty,
\qquad \mu\in(0,2],
\]
where
$\eta_{\min}:=\min\left\{\tau_{\min},\frac1{L_s+\beta}\right\},
\eta_{\max}:=\frac1\beta$ are independent of $k$. In particular, these bounds hold for every
accepted dogleg parameter $\mu^k$.
\end{lemma}

\begin{proof}
By Lemma~\ref{lem:mono} and the identity
$\eta_\mu=\bar\tau_k$ for $\mu\in(0,1]$,
we have
$\bar\tau_k\le\eta_\mu
\le\left.\eta_\mu\right|_{\mu=2}$
for all $\mu\in(0,2]$.
Lemma~\ref{lem:ssbb} gives
$\bar\tau_k\ge\tau_{\min}
\ge\min\{\tau_{\min},1/(L_s+\beta)\}=\eta_{\min}$.
At $\mu=2$, we have
$d^{\mathrm{hy}}(2)=-\mathbf B_\beta^{-1}g_k$,
and hence
\[
\left.\eta_\mu\right|_{\mu=2}
=
\frac{g_k^\top\mathbf B_\beta^{-2}g_k}
{g_k^\top\mathbf B_\beta^{-1}g_k}
\le\frac{1}{\lambda_{\min}(\mathbf B_\beta)}
\le\frac1\beta
=\eta_{\max}.
\]
The inequality follows from
$\mathbf B_\beta^{-2}
\preceq\lambda_{\min}(\mathbf B_\beta)^{-1}
\mathbf B_\beta^{-1}$
and $g_k^\top\mathbf B_\beta^{-1}g_k>0$.
Combining these estimates proves the result.
\end{proof}

\begin{remark}
The bounds in Lemma~\ref{lem:eta} depend only
on $\tau_{\min}$, $L_s$, and $\beta$.
They apply to the local steplengths associated
with the dogleg path.
If the dogleg search is skipped or fails,
the stored steplength is set to the fixed
safeguard steplength $\eta$, which need not
belong to $[\eta_{\min},\eta_{\max}]$.
Thus, the stored steplength always belongs to
$[\min\{\eta,\eta_{\min}\},
\max\{\eta,\eta_{\max}\}]$.
\end{remark}

% =====================================================================
We first establish majorization along the dogleg
line for the Euclidean kernel $\hat h_0$.
Proposition~\ref{lem:init} then provides a sufficient
condition for the corresponding property under
a general Bregman kernel.
Fix an iteration $k$ with $g_k\neq0$ and
a parameter $\mu\in(0,2]$, and let
$d:=d^{\mathrm{hy}}(\mu)$.
Using the original coordinates, consider the line
\[
\mathcal X_\mu
:=\{y^k+t d:t\in\mathbb R\}.
\]
Since $s$ is quadratic, the change in its value along this
line is
\begin{equation}
S_\mu(t)
:=s(y^k+t d)-s(y^k)
=t\langle g_k,d\rangle
+\frac{t^2}{2}d^\top\mathbf M d.
\label{eq:quadratic_along_dogleg}
\end{equation}
For the Euclidean kernel
$\hat h_0(x)=\frac12\|x\|^2$, define
\[
T_\mu(t)
:=m_\mu^{\hat h_0}(y^k+t d;y^k)-s(y^k)
=t\langle g_\mu,d\rangle
+\frac{t^2}{2\eta_\mu}\|d\|^2.
\]
Majorization along $\mathcal X_\mu$ therefore
reduces to comparing the scalar quadratic
functions $T_\mu$ and $S_\mu$.

\begin{proposition}[Bregman majorization along the dogleg line]
\label{lem:init}
Suppose that Assumptions~\ref{ass:main}
and~\ref{ass:breg} hold.
Fix an iteration $k$ with $g_k\ne0$ and
$y^k\in\operatorname{int}\operatorname{dom}\hat h$,
and a parameter $\mu\in(0,2]$.
Define
$S_\mu(t)
:=s\bigl(y^k+t d^{\mathrm{hy}}(\mu)\bigr)-s(y^k),$ and
$T_\mu^{\hat h}(t)
:=m_\mu^{\hat h}
\bigl(y^k+t d^{\mathrm{hy}}(\mu);y^k\bigr)-s(y^k).$ Then:
\begin{enumerate}
\item[(i)]
The two functions agree to first order at $t=0$:
$T_\mu^{\hat h}(0)=S_\mu(0)=0$ and
$(T_\mu^{\hat h})'(0)=S_\mu'(0)
=\langle g_k,d^{\mathrm{hy}}(\mu)\rangle$.

\item[(ii)]
$T_\mu^{\hat h}$ is strictly convex on
$\{t\in\mathbb R:
y^k+t d^{\mathrm{hy}}(\mu)
\in\operatorname{int}\operatorname{dom}\hat h\}$.

\item[(iii)]
If
\[
\frac{\lambda_{\min}^{\hat h}}{\eta_\mu}
\|d^{\mathrm{hy}}(\mu)\|^2
\ge
d^{\mathrm{hy}}(\mu)^\top\mathbf M d^{\mathrm{hy}}(\mu),
\]
then $T_\mu^{\hat h}(t)\ge S_\mu(t)$
for every $t$ in this set.
In particular, $\eta_\mu\le\lambda_{\min}^{\hat h}/L_s$
is sufficient.
\end{enumerate}
\end{proposition}

\begin{proof}
By Lemma~\ref{lem:descent},
$\langle g_k,d^{\mathrm{hy}}(\mu)\rangle<0$.
Thus, $d^{\mathrm{hy}}(\mu)\ne0$ and $\eta_\mu>0$.

\noindent\textbf{(i).}
By the definition of the local Bregman model
in \eqref{eq:dogleg_local_model}
and the definition of $T_\mu^{\hat h}$, we have
\[
T_\mu^{\hat h}(t)
=
t\langle g_\mu,d^{\mathrm{hy}}(\mu)\rangle
+\frac{1}{\eta_\mu}
D_{\hat h}\bigl(
y^k+t d^{\mathrm{hy}}(\mu),y^k
\bigr).
\]
The Bregman term and its derivative with respect
to $t$ both vanish at $t=0$.
Together with \eqref{eq:projected_slope_identity},
this gives
$T_\mu^{\hat h}(0)=S_\mu(0)=0$ and
$(T_\mu^{\hat h})'(0)=S_\mu'(0)
=\langle g_k,d^{\mathrm{hy}}(\mu)\rangle$.

\noindent\textbf{(ii).}
By Assumption~\ref{ass:breg},
\begin{align*}
(T_\mu^{\hat h})''(t)
&=
\frac{
d^{\mathrm{hy}}(\mu)^\top
\nabla^2\hat h\bigl(
y^k+t d^{\mathrm{hy}}(\mu)
\bigr)d^{\mathrm{hy}}(\mu)
}{\eta_\mu}\ge
\frac{\lambda_{\min}^{\hat h}}{\eta_\mu}
\|d^{\mathrm{hy}}(\mu)\|^2>0.
\end{align*}
The stated set is an interval because
$\operatorname{int}\operatorname{dom}\hat h$
is convex.
Therefore, $T_\mu^{\hat h}$ is strictly convex
on this set.

\noindent\textbf{(iii).}
The quadratic structure of $s$ gives
$S_\mu(t)
=t\langle g_k,d^{\mathrm{hy}}(\mu)\rangle
+\frac{t^2}{2}
d^{\mathrm{hy}}(\mu)^\top
\mathbf M d^{\mathrm{hy}}(\mu)$.
By \eqref{eq:projected_slope_identity},
\begin{equation}\label{eq:bregman_line_difference}
\begin{aligned}
T_\mu^{\hat h}(t)-S_\mu(t)
&=
\frac{1}{\eta_\mu}
D_{\hat h}\bigl(
y^k+t d^{\mathrm{hy}}(\mu),y^k
\bigr)-
\frac{t^2}{2}
d^{\mathrm{hy}}(\mu)^\top
\mathbf M d^{\mathrm{hy}}(\mu).
\end{aligned}
\end{equation}
For every $t$ in the stated set,
the segment joining $y^k$ and
$y^k+t d^{\mathrm{hy}}(\mu)$ lies in
$\operatorname{int}\operatorname{dom}\hat h$.
The integral form of the Taylor remainder
and Assumption~\ref{ass:breg} therefore give
\begin{align*}
D_{\hat h}\bigl(
y^k+t d^{\mathrm{hy}}(\mu),y^k
\bigr)
&=
t^2\int_0^1(1-\theta)\,
d^{\mathrm{hy}}(\mu)^\top
\nabla^2\hat h\bigl(
y^k+\theta t d^{\mathrm{hy}}(\mu)
\bigr)d^{\mathrm{hy}}(\mu)\,d\theta
\ge
\frac{\lambda_{\min}^{\hat h}}{2}
t^2\|d^{\mathrm{hy}}(\mu)\|^2.
\end{align*}
This estimate holds for either sign of $t$.
Combining it with the bound in part~(iii)
in \eqref{eq:bregman_line_difference}, we obtain
\[
T_\mu^{\hat h}(t)-S_\mu(t)
\ge
\frac{t^2}{2}
\left(
\frac{\lambda_{\min}^{\hat h}}{\eta_\mu}
\|d^{\mathrm{hy}}(\mu)\|^2
-d^{\mathrm{hy}}(\mu)^\top
\mathbf M d^{\mathrm{hy}}(\mu)
\right)\ge0,
\]
where the last inequality is the condition in part~(iii).
This completes the proof.
\end{proof}

For a full-domain kernel, including the diagonal
quadratic kernel $h_{\mathbf M_{\mathrm d}}$
used in Algorithm~\ref{Algorithm3},
the segment domain requirement in
Proposition~\ref{lem:init} is automatic.
Under the condition in part~(iii)
of that proposition \ref{lem:init},
majorization holds
along the dogleg line.
The proposition concerns points on this line; the
actual proximal trial point may not lie on it and is
checked directly by \eqref{eq:dogleg_acceptance}.
The regularized-Newton dogleg candidate is used
only if it satisfies both acceptance conditions in
\eqref{eq:dogleg_acceptance}. Otherwise, the
safeguarded Bregman proximal candidate is selected.

We next establish sufficient decrease of the associated Lyapunov function.

\subsection{Lyapunov Descent and Criticality of Cluster Points}
\label{sec:lyap}

Let $\lambda_{\min}^{\hat h}$ and
$\lambda_{\max}^{\hat h}$ be the constants
in Assumption~\ref{ass:breg}.
Let $\gamma\in(0,1)$, and let
$\eta_{\min}$ and $\eta_{\max}$
be the uniform local steplength bounds
in Lemma~\ref{lem:eta}.
Observe that $\mathbf B_\beta$ enters the estimates
below only through the bounds on $\eta_\mu$;
all curvature terms use the true Hessian $\mathbf M$.
For the fixed safeguard steplength
$0<\eta<\lambda_{\min}^{\hat h}/L_s$, define
\[
\mathcal A_{\mathrm d}
:=\frac{(1-\gamma)\lambda_{\min}^{\hat h}}
{2\gamma\eta_{\max}},
\quad
\mathcal B_{\mathrm d}
:=\frac{\lambda_{\max}^{\hat h}}
{2\gamma\eta_{\min}},
\quad
\mathcal A_{\mathrm v}
:=\frac12\left(
\frac{\lambda_{\min}^{\hat h}}{\eta}-L_s
\right),
\quad
\mathcal B_{\mathrm v}
:=\frac12\left(
\frac{\lambda_{\max}^{\hat h}}{\eta}
-\lambda_{\min}(\mathbf M)
\right).
\]
All four constants are positive.
Indeed, since $\eta_{\max}=1/\beta$,
\[
\mathcal A_{\mathrm d}
=\frac{(1-\gamma)\lambda_{\min}^{\hat h}\beta}
{2\gamma}>0,
\]
and $\mathcal B_{\mathrm d}>0$ because
$\gamma\in(0,1)$ and $\eta_{\min}>0$.
The restriction on $\eta$ gives
$\mathcal A_{\mathrm v}>0$ and
$\lambda_{\max}^{\hat h}/\eta
\ge\lambda_{\min}^{\hat h}/\eta
>L_s\ge\lambda_{\min}(\mathbf M)$,
which also implies $\mathcal B_{\mathrm v}>0$;
here $\lambda_{\min}(\mathbf M)$ may be zero.
Set
\begin{equation}\label{eq:AB-uniform}
\mathcal A
:=\min\{\mathcal A_{\mathrm d},\mathcal A_{\mathrm v}\}>0,
\qquad
\mathcal B
:=\max\{\mathcal B_{\mathrm d},\mathcal B_{\mathrm v}\}
<+\infty.
\end{equation}
These constants are uniform in $k$.

Recall that $\Delta_k:=x^k-x^{k-1}$,
with $x^{-1}=x^0$.
For the extrapolation parameter $\zeta\in[0,1)$,
let $\delta:=\mathcal A\zeta(1-\zeta)
+\mathcal B\zeta^2$ and define
\begin{equation}\label{eq:lifted-H}
\mathcal H(x,z):=Q(x)+\delta\|x-z\|^2,
\qquad
H_k:=\mathcal H(x^k,x^{k-1}).
\end{equation}

\begin{lemma}[Uniform Lyapunov descent]\label{lem:lyap}
Suppose that Assumptions~\ref{ass:main}
and~\ref{ass:breg} hold.
Let $\{x^k\}$ be generated by
Algorithm~\ref{Algorithm3}
under the standing initialization and domain
assumptions, with the stopping criteria disabled.
If
\begin{equation}\label{eq:zeta-new}
0\le\zeta<
\frac{\sqrt{\mathcal A}}
{\sqrt{\mathcal A}+\sqrt{\mathcal B}},
\end{equation}
then, for every $k\ge0$,
\begin{equation}\label{eq:H-descent-new}
H_{k+1}-H_k
\le-a\|x^{k+1}-x^k\|^2,
\end{equation}
where
$a:=\mathcal A(1-\zeta)^2-\mathcal B\zeta^2>0$.
Consequently, $\{H_k\}$ is nonincreasing
and converges to a finite limit.
Moreover,
\begin{equation}\label{eq:square-summable-new}
\sum_{k=0}^{\infty}\|x^{k+1}-x^k\|^2<+\infty,
\qquad
\|x^{k+1}-x^k\|\to0,
\qquad
\|x^{k+1}-y^k\|\to0.
\end{equation}
\end{lemma}

\begin{proof}
We work in the original coordinates and estimate the dogleg and safeguard candidates separately. We first consider an admissible dogleg candidate
$\tilde x^k$, when available.
By its global minimality in
\eqref{eq:dogleg_bregman_candidate},
comparison with $x^k$ gives
\begin{equation}\label{eq:r-dogleg}
r(\tilde x^k)-r(x^k)
\le\langle g_{\mu^k},x^k-\tilde x^k\rangle
+\frac{D_{\hat h}(x^k,y^k)-D_{\hat h}(\tilde x^k,y^k)}
{\gamma\eta_{\mu^k}}.
\end{equation}
The accepted majorization condition gives
\[
s(\tilde x^k)
\le s(y^k)+\langle g_{\mu^k},\tilde x^k-y^k\rangle
+\frac{1}{\eta_{\mu^k}}D_{\hat h}(\tilde x^k,y^k).
\]
The convexity of the quadratic function $s$
also gives
$s(x^k)\ge s(y^k)+\langle g_k,x^k-y^k\rangle$.
Subtracting yields
\begin{equation}\label{eq:s-dogleg}
s(\tilde x^k)-s(x^k)
\le-\langle g_k,x^k-y^k\rangle
+\langle g_{\mu^k},\tilde x^k-y^k\rangle
+\frac{1}{\eta_{\mu^k}}D_{\hat h}(\tilde x^k,y^k).
\end{equation}
Adding \eqref{eq:r-dogleg} and \eqref{eq:s-dogleg}
leaves the linear term
$\langle g_{\mu^k}-g_k,x^k-y^k\rangle$,
which is nonpositive by the acceptance condition.
Therefore,
\[
Q(\tilde x^k)-Q(x^k)
\le-\frac{1-\gamma}{\gamma\eta_{\mu^k}}
D_{\hat h}(\tilde x^k,y^k)
+\frac{1}{\gamma\eta_{\mu^k}}D_{\hat h}(x^k,y^k).
\]
Applying the Bregman bounds and
$\eta_{\min}\le\eta_{\mu^k}\le\eta_{\max}$ gives
\[
Q(\tilde x^k)-Q(x^k)
\le-\mathcal A_{\mathrm d}\|\tilde x^k-y^k\|^2
+\mathcal B_{\mathrm d}\|x^k-y^k\|^2.
\]

\noindent\emph{Safeguard candidate.}
Global minimality in
\eqref{eq:safeguard_bregman_candidate}
gives
\begin{equation}\label{eq:r-safe}
r(v^{k+1})-r(x^k)
\le-\langle g_k,v^{k+1}-x^k\rangle
+\frac{D_{\hat h}(x^k,y^k)-D_{\hat h}(v^{k+1},y^k)}{\eta}.
\end{equation}
Since $s$ is quadratic,
\begin{equation}\label{eq:s-safe}
s(v^{k+1})-s(x^k)
=\langle g_k,v^{k+1}-x^k\rangle
+\frac12\|v^{k+1}-y^k\|_{\mathbf M}^2
-\frac12\|x^k-y^k\|_{\mathbf M}^2.
\end{equation}
Adding \eqref{eq:r-safe} and \eqref{eq:s-safe}
cancels the linear terms and yields
\[
\begin{aligned}
Q(v^{k+1})-Q(x^k)
&\le-\frac1\eta D_{\hat h}(v^{k+1},y^k)
+\frac12\|v^{k+1}-y^k\|_{\mathbf M}^2
+\frac1\eta D_{\hat h}(x^k,y^k)
-\frac12\|x^k-y^k\|_{\mathbf M}^2.
\end{aligned}
\]
Using the Bregman bounds and
$\lambda_{\min}(\mathbf M)\|z\|^2
\le\|z\|_{\mathbf M}^2\le L_s\|z\|^2$,
we obtain
\begin{equation}\label{eq:safe-uniform}
Q(v^{k+1})-Q(x^k)
\le-\mathcal A_{\mathrm v}\|v^{k+1}-y^k\|^2
+\mathcal B_{\mathrm v}\|x^k-y^k\|^2.
\end{equation}
This estimate holds whether or not the safeguard
candidate is selected.
It will also be used in Lemma~\ref{lem:safeguard}.

The next iterate is either an admissible dogleg
candidate or the safeguard candidate.
By \eqref{eq:AB-uniform}, the corresponding
estimate implies
\begin{equation}\label{eq:common-descent}
Q(x^{k+1})-Q(x^k)
\le-\mathcal A\|x^{k+1}-y^k\|^2
+\mathcal B\|x^k-y^k\|^2.
\end{equation}
The extrapolation formula gives
$x^k-y^k=-\zeta\Delta_k$ and
$x^{k+1}-y^k=\Delta_{k+1}-\zeta\Delta_k$.
Young's inequality implies
\[
\|\Delta_{k+1}-\zeta\Delta_k\|^2
\ge(1-\zeta)\|\Delta_{k+1}\|^2
-\zeta(1-\zeta)\|\Delta_k\|^2.
\]
Substituting into \eqref{eq:common-descent}
and using the definition of $\delta$ gives
\[
Q(x^{k+1})-Q(x^k)
\le-\mathcal A(1-\zeta)\|\Delta_{k+1}\|^2
+\delta\|\Delta_k\|^2.
\]
Hence, $H_{k+1}-H_k
\le-\bigl[\mathcal A(1-\zeta)-\delta\bigr]
\|\Delta_{k+1}\|^2
=-a\|\Delta_{k+1}\|^2.$
Condition~\eqref{eq:zeta-new} ensures $a>0$.
Since $\delta\ge0$ and $Q(x)\ge\underline Q$,
we have $H_k\ge\underline Q$.
Thus, $\{H_k\}$ is nonincreasing and convergent.
Summation gives
\[
a\sum_{k=0}^{N}\|\Delta_{k+1}\|^2
\le H_0-H_{N+1}
\le H_0-\underline Q.
\]
Letting $N\to\infty$ proves square summability,
which implies $\Delta_k\to0$.
Finally,
$\|x^{k+1}-y^k\|
\le\|\Delta_{k+1}\|+\zeta\|\Delta_k\|\to0$.
The proof uses global minimality of the proximal
candidates and does not require convexity of $r$.
\end{proof}

\begin{remark}[Instance conditions for the descent lemma]
\label{rem:parameter-check}
Lemma~\ref{lem:lyap} requires
\begin{equation}\label{eq:zeta-check}
\zeta<
\frac{\sqrt{\mathcal A}}{\sqrt{\mathcal A}+\sqrt{\mathcal B}}.
\end{equation}
Two conditions must therefore be checked for every
experimental instance:
\begin{enumerate}
\item[(P1)] $\eta<\lambda_{\min}^{\hat h}/L_s$;
for the choice $\eta=0.99/L_s$ this is
$\lambda_{\min}^{\hat h}>0.99$;
\item[(P2)] $\zeta$ is below the threshold in
\eqref{eq:zeta-check}.
\end{enumerate}
Condition (P2) is sensitive to $\beta$.
Since $\eta_{\max}=1/\beta$,
\[
\mathcal A_{\mathrm d}
=\frac{(1-\gamma)\lambda_{\min}^{\hat h}\beta}
{2\gamma},
\]
so $\mathcal A\le\mathcal A_{\mathrm d}=O(\beta)$.
For $\beta=10^{-6}$ and $\gamma=0.998$,
$\sqrt{\mathcal A}\lesssim 3\times10^{-5}$.
On the other hand,
\[
\mathcal B_{\mathrm v}
=\frac12\left(
\frac{\lambda_{\max}^{\hat h}}{\eta}
-\lambda_{\min}(\mathbf M)
\right)
\approx\frac{\lambda_{\max}^{\hat h}L_s}{1.98},
\]
so $\sqrt{\mathcal B}$ is of order one for the
instances considered. Hence the threshold in
\eqref{eq:zeta-check} is of order $10^{-5}$,
whereas the experiments use $\zeta=0.02$ and
$\zeta=0.045$. These choices are far above the
sufficient bound; the experiments are therefore
not covered by Lemma~\ref{lem:lyap} unless they
are rerun with a smaller $\zeta$ or a larger
$\beta$. Otherwise they should be presented as
empirical results.
\end{remark}

\begin{lemma}[Residual bound and asymptotic proximity
of the safeguard point]\label{lem:safeguard}
Under the assumptions of Lemma~\ref{lem:lyap},
the safeguard candidate satisfies
\begin{equation}\label{eq:safe-stationarity}
\operatorname{dist}\bigl(0,\partial Q(v^{k+1})\bigr)
\le
\left(L_s+\frac{\lambda_{\max}^{\hat h}}{\eta}\right)
\|v^{k+1}-y^k\|.
\end{equation}
Moreover,
$\|v^{k+1}-y^k\|\to0$,
$\|v^{k+1}-x^k\|\to0$, and
$\operatorname{dist}(0,\partial Q(v^{k+1}))\to0$.
\end{lemma}

\begin{proof}
The optimality condition of the safeguard
subproblem gives
\[
w^{k+1}
:=\nabla s(v^{k+1})-g_k
-\frac1\eta
\bigl[
\nabla\hat h(v^{k+1})-\nabla\hat h(y^k)
\bigr]
\in\partial Q(v^{k+1}).
\]
By Assumption~\ref{ass:breg} and the convexity
of $\operatorname{int}\operatorname{dom}\hat h$,
the Hessian upper bound gives
$\|\nabla\hat h(v^{k+1})-\nabla\hat h(y^k)\|
\le\lambda_{\max}^{\hat h}\|v^{k+1}-y^k\|$.
Since $g_k=\nabla s(y^k)$ and $\nabla s$
is $L_s$-Lipschitz continuous, we obtain
\[
\operatorname{dist}\bigl(0,\partial Q(v^{k+1})\bigr)
\le\|w^{k+1}\|
\le
\left(L_s+\frac{\lambda_{\max}^{\hat h}}{\eta}\right)
\|v^{k+1}-y^k\|,
\]
which proves \eqref{eq:safe-stationarity}.

Estimate~\eqref{eq:safe-uniform} holds whether
or not the safeguard candidate is selected.
Combining it with the selection rule
$Q(x^{k+1})\le Q(v^{k+1})$ yields
\begin{equation}\label{eq:safe-selection-bound}
\mathcal A_{\mathrm v}\|v^{k+1}-y^k\|^2
\le Q(x^k)-Q(x^{k+1})
+\mathcal B_{\mathrm v}\zeta^2\|\Delta_k\|^2.
\end{equation}
By Lemma~\ref{lem:lyap}, $H_k$ converges
and $\Delta_k\to0$.
Since $Q(x^k)=H_k-\delta\|\Delta_k\|^2$,
the sequence $\{Q(x^k)\}$ also converges.
Hence, $Q(x^k)-Q(x^{k+1})\to0$,
and the right-hand side of
\eqref{eq:safe-selection-bound} tends to zero.
As $\mathcal A_{\mathrm v}>0$,
it follows that $\|v^{k+1}-y^k\|\to0$.
Moreover,
$\|v^{k+1}-x^k\|
\le\|v^{k+1}-y^k\|+\zeta\|\Delta_k\|\to0$.
Finally, \eqref{eq:safe-stationarity} implies
$\operatorname{dist}(0,\partial Q(v^{k+1}))\to0$.
\end{proof}

\begin{theorem}[Criticality of cluster points]
\label{thm:criticality}
Under the assumptions of Lemma~\ref{lem:lyap},
every cluster point $x^\star$ of $\{x^k\}$
belonging to
$\operatorname{int}\operatorname{dom}\hat h$
is critical for $Q$, that is,
$0\in\partial Q(x^\star)$.
In particular, for a full-domain kernel,
every cluster point of $\{x^k\}$ is critical.
\end{theorem}

\begin{proof}
Let $x^{k_j}\to x^\star$, where
$x^\star\in\operatorname{int}\operatorname{dom}\hat h$.
Lemmas~\ref{lem:lyap} and~\ref{lem:safeguard}
give $y^{k_j}\to x^\star$ and
$v^{k_j+1}\to x^\star$.
The subgradients constructed in the proof of
Lemma~\ref{lem:safeguard} satisfy
$w^{k_j+1}\in\partial Q(v^{k_j+1})$
and $w^{k_j+1}\to0$.
To apply the closedness property of the limiting
subdifferential, it remains to establish
$Q(v^{k_j+1})\to Q(x^\star)$.

Lower semicontinuity of $Q$ and
$Q(x^{k_j})\le H_{k_j}\le H_0<+\infty$
imply
$Q(x^\star)\le
\liminf_{j\to\infty}Q(x^{k_j})
\le H_0<+\infty$.
Thus, $x^\star\in\operatorname{dom}r$.
By global minimality of the safeguard candidate,
comparison with $x^\star$ gives
\[
r(v^{k_j+1})-r(x^\star)
\le
\langle g_{k_j},x^\star-v^{k_j+1}\rangle
+\frac{
D_{\hat h}(x^\star,y^{k_j})
-D_{\hat h}(v^{k_j+1},y^{k_j})
}{\eta}.
\]
Since $g_{k_j}=\nabla s(y^{k_j})
\to\nabla s(x^\star)$ and
$v^{k_j+1}\to x^\star$,
the inner product tends to zero.
By Assumption~\ref{ass:breg},
\[
0\le D_{\hat h}(x^\star,y^{k_j})
\le\frac{\lambda_{\max}^{\hat h}}{2}
\|x^\star-y^{k_j}\|^2\to0.
\]
Using also
$D_{\hat h}(v^{k_j+1},y^{k_j})\ge0$,
we obtain
$\limsup_{j\to\infty}r(v^{k_j+1})
\le r(x^\star)$.
Lower semicontinuity of $r$ gives
$r(x^\star)\le\liminf_{j\to\infty}r(v^{k_j+1})$,
and hence $r(v^{k_j+1})\to r(x^\star)$.
By continuity of $s$,
$Q(v^{k_j+1})\to Q(x^\star)$.

Therefore,
$v^{k_j+1}\to x^\star$,
$Q(v^{k_j+1})\to Q(x^\star)$, and
$w^{k_j+1}\in\partial Q(v^{k_j+1})$
with $w^{k_j+1}\to0$.
The closedness of the limiting subdifferential
along function-value convergent sequences
then yields $0\in\partial Q(x^\star)$.
\end{proof}

\subsection{K{\L}-Based Convergence and Finite Length}
\label{sec:kl}

\begin{lemma}[Relative error bounds]\label{lem:relerr}
Under the assumptions of Lemma~\ref{lem:lyap}, let
$\bar C:=L_s+
\frac{\lambda_{\max}^{\hat h}}
{\min\{\gamma\eta_{\min},\eta\}}$.
Define $e^k:=g_{\mu^k}-g_k$ when the dogleg
candidate is selected as $x^{k+1}$,
and $e^k:=0$ otherwise.
Then
\begin{equation}\label{eq:Q-relative}
\operatorname{dist}\bigl(0,\partial Q(x^{k+1})\bigr)
\le\bar C\|x^{k+1}-y^k\|+\|e^k\|.
\end{equation}
For the Lyapunov function defined in
\eqref{eq:lifted-H}, we also have
\begin{equation}\label{eq:H-relative}
\operatorname{dist}\bigl(
0,\partial\mathcal H(x^{k+1},x^k)
\bigr)
\le(\bar C+4\delta)\|\Delta_{k+1}\|
+\bar C\zeta\|\Delta_k\|+\|e^k\|.
\end{equation}
\end{lemma}

\begin{proof}
If the dogleg candidate is selected,
its optimality condition gives
\[
\nabla s(x^{k+1})-g_{\mu^k}
-\frac{
\nabla\hat h(x^{k+1})-\nabla\hat h(y^k)
}{\gamma\eta_{\mu^k}}
\in\partial Q(x^{k+1}).
\]
Since $g_{\mu^k}=g_k+e^k$ and
$g_k=\nabla s(y^k)$, the Lipschitz continuity
of $\nabla s$ and the Hessian upper bound
in Assumption~\ref{ass:breg} yield
\[
\operatorname{dist}\bigl(0,\partial Q(x^{k+1})\bigr)
\le
\left(
L_s+\frac{\lambda_{\max}^{\hat h}}
{\gamma\eta_{\mu^k}}
\right)
\|x^{k+1}-y^k\|+\|e^k\|.
\]
Using $\eta_{\mu^k}\ge\eta_{\min}$ and
the definition of $\bar C$, we obtain
\eqref{eq:Q-relative}.

If the safeguard candidate is selected,
its optimality condition similarly gives
\[
\operatorname{dist}\bigl(0,\partial Q(x^{k+1})\bigr)
\le
\left(
L_s+\frac{\lambda_{\max}^{\hat h}}{\eta}
\right)
\|x^{k+1}-y^k\|.
\]
Since $e^k=0$ in this case,
\eqref{eq:Q-relative} holds in both cases. By the smooth perturbation rule for the
limiting subdifferential,
$\partial\mathcal H(x^{k+1},x^k)
=
\left\{
\bigl(
\xi+2\delta\Delta_{k+1},
-2\delta\Delta_{k+1}
\bigr):
\xi\in\partial Q(x^{k+1})
\right\}.$ For any $\xi\in\partial Q(x^{k+1})$,
the triangle inequality in the Euclidean
product space and $\delta\ge0$ give
\begin{align*}
\bigl\|
(\xi+2\delta\Delta_{k+1},
-2\delta\Delta_{k+1})
\bigr\|
&\le
\|\xi+2\delta\Delta_{k+1}\|
+2\delta\|\Delta_{k+1}\|\\
&\le
\|\xi\|+4\delta\|\Delta_{k+1}\|.
\end{align*}
Taking the infimum over $\xi$ yields
$\operatorname{dist}\bigl(
0,\partial\mathcal H(x^{k+1},x^k)
\bigr)
\le
\operatorname{dist}\bigl(0,\partial Q(x^{k+1})\bigr)
+4\delta\|\Delta_{k+1}\|.$  Finally, the extrapolation formula implies
$\|x^{k+1}-y^k\|
\le\|\Delta_{k+1}\|+\zeta\|\Delta_k\|$.
Combining this inequality with
\eqref{eq:Q-relative} proves
\eqref{eq:H-relative}.
\end{proof}
For the finite-length and convergence analysis, we make the following
additional assumptions: \textup{(C1)} the sequence $\{x^k\}$ is
bounded; \textup{(C2)} the Lyapunov function $\mathcal H$ is a
K{\L} function; \textup{(C3)} the gradient mismatch sequence
defined in Lemma~\ref{lem:relerr} satisfies
$\sum_{k=0}^{\infty}\|e^k\|<+\infty$; and \textup{(C4)} there
exists $\varkappa>0$ such that
$\|\Delta_k\|\le\varkappa\|\Delta_{k+1}\|$ for all sufficiently
large $k$.

Substituting the bound in \textup{(C4)} into \eqref{eq:H-relative}
gives, for all sufficiently large $k$,
\[
\operatorname{dist}\bigl(
0,\partial\mathcal H(x^{k+1},x^k)
\bigr)
\le
\bigl[\bar C(1+\varkappa\zeta)+4\delta\bigr]
\|\Delta_{k+1}\|+\|e^k\|.
\]

\begin{theorem}[Conditional finite length]\label{thm:finlen}
Under the assumptions of Lemma~\ref{lem:lyap},
suppose that (C1)--(C4) hold and that every cluster
point of $\{x^k\}$ belongs to
$\operatorname{int}\operatorname{dom}\hat h$.
Then $\{x^k\}$ has finite length,
$\sum_{k=0}^{\infty}\|x^{k+1}-x^k\|<+\infty$,
and converges to a critical point of $Q$.
\end{theorem}

\begin{proof}
Set $z^k:=(x^k,x^{k-1})$, so that
$\mathcal H(z^k)=H_k$.
For all sufficiently large $k$, the ratio
condition gives
\[
\|z^{k+1}-z^k\|^2
=\|\Delta_{k+1}\|^2+\|\Delta_k\|^2
\le(1+\varkappa^2)\|\Delta_{k+1}\|^2.
\]
Consequently, Lemma~\ref{lem:lyap} implies
\begin{equation}\label{eq:lifted-descent}
\mathcal H(z^k)-\mathcal H(z^{k+1})
\ge\frac{a}{1+\varkappa^2}
\|z^{k+1}-z^k\|^2,
\end{equation}
where $a>0$ is the constant in
\eqref{eq:H-descent-new}. By Lemma~\ref{lem:relerr} and the ratio condition,
\[
\operatorname{dist}\bigl(
0,\partial\mathcal H(z^{k+1})
\bigr)
\le B\|z^{k+1}-z^k\|+\|e^k\|,
\qquad
B:=\bar C(1+\varkappa\zeta)+4\delta>0.
\]
Thus,
$\frac1B\operatorname{dist}\bigl(
0,\partial\mathcal H(z^{k+1})
\bigr)
\le\|z^{k+1}-z^k\|+\frac{\|e^k\|}{B}.$  Together with \eqref{eq:lifted-descent},
this verifies the sufficient-decrease and
relative-error conditions of
\cite[Theorem~1]{frankel2015splitting}
on a tail of the sequence.
The associated coefficients
$a/(1+\varkappa^2)$ and $1/B$ are fixed
and positive.
Hence, the parameter conditions of that theorem
hold, including summability of the error sequence
$\{\|e^k\|/B\}$.

It remains to verify $\mathcal H$-precompactness.
By boundedness, choose a subsequence
$x^{k_j}\to x^\star$.
By assumption,
$x^\star\in\operatorname{int}\operatorname{dom}\hat h$.
Lemma~\ref{lem:lyap} gives $\Delta_k\to0$,
so $z^{k_j}\to(x^\star,x^\star)$ and
$x^{k_j+1}\to x^\star$.
Since $Q(x^k)=H_k-\delta\|\Delta_k\|^2$,
the sequence $\{Q(x^k)\}$ converges to
$H_\infty:=\lim\limits_{k\to\infty}H_k$.
The proof of Theorem~\ref{thm:criticality}
gives $Q(v^{k_j+1})\to Q(x^\star)$.
Lower semicontinuity of $Q$ and the selection rule
$Q(x^{k_j+1})\le Q(v^{k_j+1})$
therefore yield
\[
Q(x^\star)
\le\liminf_{j\to\infty}Q(x^{k_j})
=H_\infty
=\lim_{j\to\infty}Q(x^{k_j+1})
\le\lim_{j\to\infty}Q(v^{k_j+1})
=Q(x^\star).
\]
Thus, $H_\infty=Q(x^\star)$.
Since $\mathcal H(x^\star,x^\star)=Q(x^\star)$,
we obtain
$\mathcal H(z^{k_j})
\to\mathcal H(x^\star,x^\star)$.
This establishes $\mathcal H$-precompactness.

Since $\mathcal H$ is a K{\L} function,
\cite[Theorem~1]{frankel2015splitting},
applied to a tail of the sequence,
implies that $\{z^k\}$ converges
and that its tail has finite length.
Adding the finitely many initial terms gives
\[
\sum_{k=0}^{\infty}\|x^{k+1}-x^k\|
\le
\sum_{k=0}^{\infty}\|z^{k+1}-z^k\|
<+\infty.
\]
Therefore, $\{x^k\}$ converges.
Its limit belongs to
$\operatorname{int}\operatorname{dom}\hat h$
by assumption and is critical for $Q$
by Theorem~\ref{thm:criticality}.
\end{proof}

\begin{theorem}[Conditional rate of convergence]\label{thm:kl-rates}
Suppose that the assumptions of
Theorem~\ref{thm:finlen} hold, and let
$x^\star$ denote the limit of $\{x^k\}$.
Assume that $\mathcal H$ satisfies the
K{\L} inequality at
$(x^\star,x^\star)$ with desingularizing function
$\varphi(t)=c_\varphi t^{1-\theta}$,
where $c_\varphi>0$ and $\theta\in[0,1)$.
Suppose further that there exists $c_e\ge0$
such that
\begin{equation}\label{eq:error-step-control}
\|e^k\|\le c_e\|\Delta_{k+1}\|
\end{equation}
for all sufficiently large $k$.
This condition is not implied by (C1)--(C4);
it is a separate control of the gradient mismatch
that is needed only for the rate.
Then:
\begin{enumerate}
\item[(i)]
If $\theta=0$, the sequence reaches $x^\star$
in finitely many iterations.

\item[(ii)]
If $\theta\in(0,1/2]$, there exist
$c_1>0$, $q\in(0,1)$, and $k_0\in\mathbb N$
such that
$\|x^k-x^\star\|\le c_1q^{k-k_0}$
for all $k\ge k_0$.

\item[(iii)]
If $\theta\in(1/2,1)$, there exist
$c_2>0$ and $k_0\in\mathbb N$ such that
$\|x^k-x^\star\|
\le c_2(k-k_0+1)^{-\frac{1-\theta}{2\theta-1}},  k\ge k_0$.
\end{enumerate}
\end{theorem}

\begin{proof}
Set $z^k:=(x^k,x^{k-1})$,
$z^\star:=(x^\star,x^\star)$, and
$H_\infty:=\mathcal H(z^\star)$.
By Theorem~\ref{thm:finlen} and its proof,
$z^k\to z^\star$ and $H_k\to H_\infty$.
Let $E_k:=H_k-H_\infty\ge0$.
Lemma~\ref{lem:lyap} gives
\begin{equation}\label{eq:rate-descent}
E_k-E_{k+1}
\ge a\|\Delta_{k+1}\|^2.
\end{equation}
If $E_k=0$ for some $k$, monotonicity and
\eqref{eq:rate-descent} imply that all subsequent
step differences vanish.
Since $x^k\to x^\star$, the sequence then
remains at $x^\star$.
It therefore suffices to consider the case
$E_k>0$ for every $k$.

By Lemma~\ref{lem:relerr}, the ratio condition,
and \eqref{eq:error-step-control}, there exists
$k_0$ such that
\begin{equation}\label{eq:rate-relative}
\operatorname{dist}\bigl(
0,\partial\mathcal H(z^{k+1})
\bigr)
\le C_r\|\Delta_{k+1}\|,
\qquad k\ge k_0,
\end{equation}
where
$C_r:=\bar C(1+\varkappa\zeta)+4\delta+c_e>0$.
Increasing $k_0$ if necessary, the KL inequality
holds at $z^k$ for every $k\ge k_0$.
Since $\varphi'(t)=c_\varphi(1-\theta)t^{-\theta}$,
it gives
$E_{k+1}^{\theta}
\le c_\varphi(1-\theta)
\operatorname{dist}\bigl(
0,\partial\mathcal H(z^{k+1})
\bigr)
\le c_\varphi(1-\theta)C_r
\|\Delta_{k+1}\|.$
Combining this estimate with
\eqref{eq:rate-descent} yields
\begin{equation}\label{eq:rate-recursion}
E_k-E_{k+1}\ge\nu E_{k+1}^{2\theta},
\qquad
\nu:=\frac{a}{
c_\varphi^2(1-\theta)^2C_r^2}>0.
\end{equation}

We next derive a distance estimate.
For $k\ge k_0+1$, concavity of $\varphi$,
the KL inequality at $z^k$, and
\eqref{eq:rate-relative} imply
\begin{align*}
\varphi(E_k)-\varphi(E_{k+1})
&\ge\varphi'(E_k)(E_k-E_{k+1})\ge
\frac{a\|\Delta_{k+1}\|^2}
{C_r\|\Delta_k\|}.
\end{align*}
Here $\|\Delta_k\|>0$, since otherwise
\eqref{eq:rate-relative} and the KL inequality
at $z^k$ would contradict $E_k>0$.
Using $2\sqrt{uv}\le u+v$, we obtain
$2\|\Delta_{k+1}\|
\le\|\Delta_k\|
+\frac{C_r}{a}
\bigl[\varphi(E_k)-\varphi(E_{k+1})\bigr].$
Summing from $k$ to $N$ and rearranging gives
\[
\sum_{j=k}^{N}\|\Delta_{j+1}\|
+\|\Delta_{N+1}\|
\le\|\Delta_k\|
+\frac{C_r}{a}
\bigl[\varphi(E_k)-\varphi(E_{N+1})\bigr].
\]
Letting $N\to\infty$ and using
\eqref{eq:rate-descent}, we conclude that
\begin{equation}\label{eq:rate-distance}
\|x^k-x^\star\|
\le\|\Delta_k\|+\frac{C_r}{a}\varphi(E_k)
\le\frac{1}{\sqrt a}E_{k-1}^{1/2}
+\frac{C_rc_\varphi}{a}E_k^{1-\theta},
\qquad k\ge k_0+1.
\end{equation}

\noindent\emph{(i).}
Suppose that $\theta=0$.
If $E_k$ remained positive,
\eqref{eq:rate-recursion} would give
$E_k-E_{k+1}\ge\nu$ for all $k\ge k_0$,
which contradicts $E_k\to0$.
Thus, $E_k$ reaches zero in finitely many
iterations, after which the sequence remains
at $x^\star$.

\noindent\emph{(ii).}
Suppose that $\theta\in(0,1/2]$.
Increasing $k_0$ if necessary, assume that
$E_k\le1$ for all $k\ge k_0$.
Then $E_{k+1}^{2\theta}\ge E_{k+1}$,
so \eqref{eq:rate-recursion} yields
$E_{k+1}\le E_k/(1+\nu)$ and hence
$E_k\le E_{k_0}(1+\nu)^{-(k-k_0)}$.
Since $1-\theta\ge1/2$,
\eqref{eq:rate-distance} gives
$\|x^k-x^\star\|\le c_1q^{k-k_0}$
for $k\ge k_0+1$, where
$q:=(1+\nu)^{-1/2}\in(0,1)$ and $c_1>0$
is sufficiently large.
Increasing $c_1$ if necessary extends the
estimate to $k=k_0$.

\noindent\emph{(iii).}
Suppose that $\theta\in(1/2,1)$ and let
$p:=2\theta-1>0$.
If $E_{k+1}\ge E_k/2$, then
\begin{align*}
E_{k+1}^{-p}-E_k^{-p}
&=p\int_{E_{k+1}}^{E_k}t^{-p-1}\,dt\ge p(E_k-E_{k+1})E_k^{-p-1}
\ge p\nu
\left(\frac{E_{k+1}}{E_k}\right)^{p+1}
\ge p\nu\,2^{-(p+1)}.
\end{align*}
Otherwise, $E_{k+1}<E_k/2$, and hence
$E_{k+1}^{-p}-E_k^{-p}
>(2^p-1)E_k^{-p}
\ge(2^p-1)E_{k_0}^{-p}.$ Thus, with
$c_0:=\min\{p\nu\,2^{-(p+1)},
(2^p-1)E_{k_0}^{-p}\}>0$,
we have
$E_{k+1}^{-p}-E_k^{-p}\ge c_0$
for all $k\ge k_0$.
Summation gives
\[
E_k
\le
\bigl[E_{k_0}^{-p}+c_0(k-k_0)\bigr]^{-1/p}
=O\left(k^{-\frac1{2\theta-1}}\right).
\]
The two terms on the right-hand side of
\eqref{eq:rate-distance} are therefore
$O(k^{-1/[2(2\theta-1)]})$ and
$O(k^{-(1-\theta)/(2\theta-1)})$, respectively.
Since $1-\theta<1/2$, this yields
$\|x^k-x^\star\|
=O(k^{-(1-\theta)/(2\theta-1)})$.
Choosing $c_2>0$ sufficiently large gives
\[
\|x^k-x^\star\|
\le c_2(k-k_0+1)^{-\frac{1-\theta}{2\theta-1}},
\qquad k\ge k_0.
\]
This completes the proof.
\end{proof}

\begin{remark}[Scope of the convergence results]
Under the assumptions of Lemma~\ref{lem:lyap},
Theorem~\ref{thm:criticality} establishes criticality of
all cluster points in the interior of the kernel domain
without requiring (C1)--(C4) or the error bound
\eqref{eq:error-step-control}.
The diagonal quadratic kernel used in
Section~\ref{sec:numerical_experiment} has domain
$\mathbb R^n$, so the interior condition holds for every
cluster point.
Moreover, the objectives in
Section~\ref{sec:numerical_experiment} are semialgebraic,
and adding the quadratic term in \eqref{eq:lifted-H}
preserves semialgebraicity.
The resulting proper, lower semicontinuous Lyapunov
function therefore satisfies the K{\L} property, which
verifies (C2) for these problems.
\end{remark}

% ---------------------------------------------------------------------

\section{Experiments}\label{sec:numerical_experiment}

We evaluate SHBB-PDOME (Algorithm~\ref{Algorithm3}) on
nonconvex sparse recovery, sparse approximation, and robust
principal component analysis (RPCA). We compare it with
PG \cite{combettes2011proximal},
mAPG \cite{li2015accelerated},
PANOC \cite{stella2017simple},
PANOC+ \cite{de2022proximal},
PDOM \cite{zhou2024proximal},
and PDOME \cite{zhao2025accelerated}.
All methods were implemented in MATLAB and run on a
Windows 10 computer equipped with an Intel Core i7-12700
processor at 2.10 GHz and 16 GB of RAM.
The sparse recovery and sparse approximation experiments
use the quadratic loss $s(x)=\frac12\|Ax-y\|^2$, with
$\nabla s(x)=A^\top(Ax-y)$,
$\mathbf M=A^\top A$, and
$L_s=\lambda_{\max}(A^\top A)$.
The sensing matrix and data generation procedures are
specified separately for each experiment.
For SHBB-PDOME, we construct the fixed diagonal metric
$\mathbf M_{\mathrm d}=\operatorname{Diag}(m_1,\ldots,m_n)$
from the sensing matrix by setting
\begin{equation}\label{eq:dct-diagonal-metric}
m_i=
\begin{cases}
\displaystyle
\left(\sum_{j=1}^{m}A_{ji}^2\right)^{\alpha_{\mathrm d}}
+\beta_{\mathrm d},
& \alpha_{\mathrm d}>0,\\[6pt]
1+\beta_{\mathrm d},
& \alpha_{\mathrm d}=0,
\end{cases}
\qquad
\alpha_{\mathrm d}\ge0,\quad \beta_{\mathrm d}>0.
\end{equation}
The metric is computed once per instance and remains fixed
throughout the run. The theoretical analysis allows any
fixed positive definite diagonal metric and does not rely
on this particular construction.
The corresponding quadratic kernel
$\hat h(x)=\frac12x^\top\mathbf M_{\mathrm d}x$
induces the Bregman distance
$D_{\hat h}(x,z)=\frac12\|x-z\|_{\mathbf M_{\mathrm d}}^2$,
with curvature bounds
$\lambda_{\min}^{\hat h}=\min_i m_i$ and
$\lambda_{\max}^{\hat h}=\max_i m_i$.

For RPCA, we use the product Frobenius norm and set
$\mathbf M_{\mathrm d}=I$.
The proximal subproblems are solved by coordinatewise
thresholding in the sparse recovery and sparse approximation
experiments, and by truncated SVD and entrywise hard
thresholding in the RPCA experiment.
The regularization parameter $\beta>0$ in
\eqref{eq:precond}, the lower steplength bound $\tau_{\min}$
in \eqref{eq:tauk}, and the remaining algorithm parameters
are specified separately for each experiment.
The convergence guarantees apply under the assumptions
of the corresponding theoretical results.

Within each trial, all methods use the same data and
initial point. The sparse recovery and sparse approximation
experiments use a common randomly generated initial point.
For RPCA, the common initialization is specified in its
experimental setup. All methods terminate when
\[
\frac{\|x^{k+1}-x^k\|}{1+\|x^{k+1}\|}
<\varepsilon_{\mathrm{tol}}
\]
or when the iteration count reaches $2000$.
We set $\varepsilon_{\mathrm{tol}}=10^{-14}$ for sparse
recovery with $\ell_0$ regularization, and
$\varepsilon_{\mathrm{tol}}=10^{-6}$ for sparse approximation
with $\ell_{1/2}$ regularization and for RPCA.
For RPCA, $x^k$ denotes the vectorized pair
$(\operatorname{vec}(L^k),\operatorname{vec}(S^k))$. For the sparse recovery and sparse approximation problems,
we report the normalized recovery error
$\mathrm{NRE}(k):=\|x^k-x^\star\|/\|x^\star\|$,
the data residual $\|Ax^k-y\|$, and iteration counts.
For RPCA, we report the relative low-rank recovery error
$\|L^k-L^\star\|_F/\|L^\star\|_F$ and iteration counts.
Recovery success criteria and the numbers of independent
trials are specified separately for each experiment.

\subsection{Nonconvex Sparse Recovery Problem} Sparse recovery here seeks to reconstruct the original sparse signal from an under-determined set of measurements.
\begin{align}
\displaystyle
\min_{x\in\mathbb{R}^n}\frac{1}{2}\|y-{\Phi}{\Upsilon}{x}\|^2+\lambda_1\|{x}\|_0,
\label{eq:44}
\end{align}
where $\Phi\in\mathbb{R}^{m\times n}$ is a row-subsampling matrix, $\Upsilon\in\mathbb{R}^{n\times n}$ is the sparse transformation basis and $\lambda_1=0.02\|A^\top y\|_\infty$. In this experiment, the sensing matrix is constructed as
$A=\Phi\Upsilon$, where
$\Phi\in\mathbb{R}^{m\times n}$ is a row-subsampling
matrix (each row has a single entry equal to one), and
$\Upsilon\in\mathbb{R}^{n\times n}$ is the discrete
cosine transform (DCT) matrix
\cite{stankovic2018analysis}. Since $\Phi$ is a subsampling matrix, $\Phi^\top\Phi$ is diagonal. The objective Hessian in \eqref{eq:44} is $\mathbf M:=\Upsilon^\top\Phi^\top\Phi\Upsilon$. It remains rank deficient. The exact regularized Hessian $\mathbf B_\beta=\mathbf M+\beta I$ is used only to construct the Newton endpoint \eqref{eq:newton-leg}; the objective is unchanged.
The penalty is separable. For a step size $a>0$,
the $i$-th coordinate of the metric-weighted proximal
mapping is
\[
x_i=
\begin{cases}
c_i, & |c_i|>\sqrt{\dfrac{2a\lambda_1}{m_i}},\\[2mm]
0, & |c_i|\le\sqrt{\dfrac{2a\lambda_1}{m_i}},
\end{cases}
\]
where $a=\gamma\eta_\mu$ for the dogleg candidate
and $a=\eta$ for the safeguard candidate.
For SHBB-PDOME, we set
$\gamma=0.998$, $\zeta=0.02$, $w=0.2$,
$\rho_0=0.5$, $\delta_\rho=0.05$,
$\alpha_{\mathrm d}=0.15$,
$\beta_{\mathrm d}=0.1$,
$\beta=0.1$, and $\eta=0.99/L_s$.
This choice of $\eta$ satisfies
$\eta<\lambda_{\min}^{\hat h}/L_s$
whenever $\lambda_{\min}^{\hat h}>0.99$.

\begin{figure}[htbp]
  \centering
  \begin{minipage}[t]{0.32\textwidth}
    \centering
    \includegraphics[width=\linewidth]{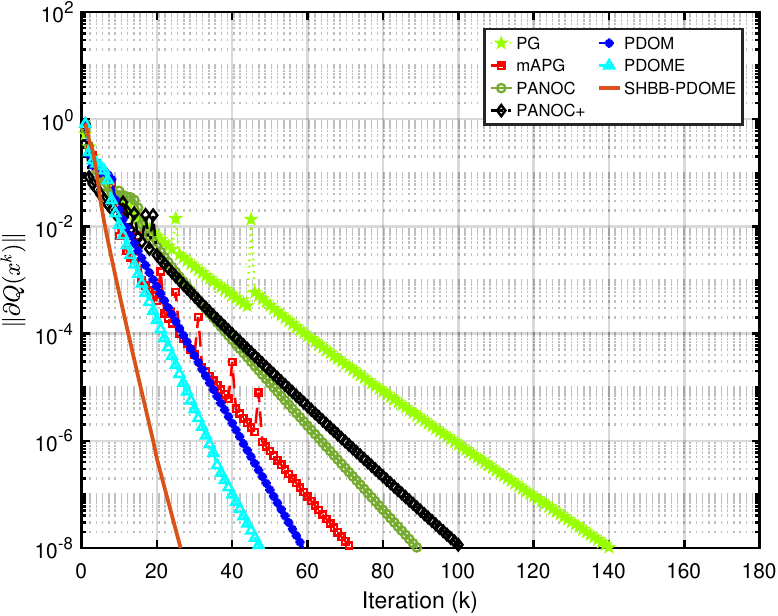}\\
    (a)
  \end{minipage}\hfill
  \begin{minipage}[t]{0.32\textwidth}
    \centering
    \includegraphics[width=\linewidth]{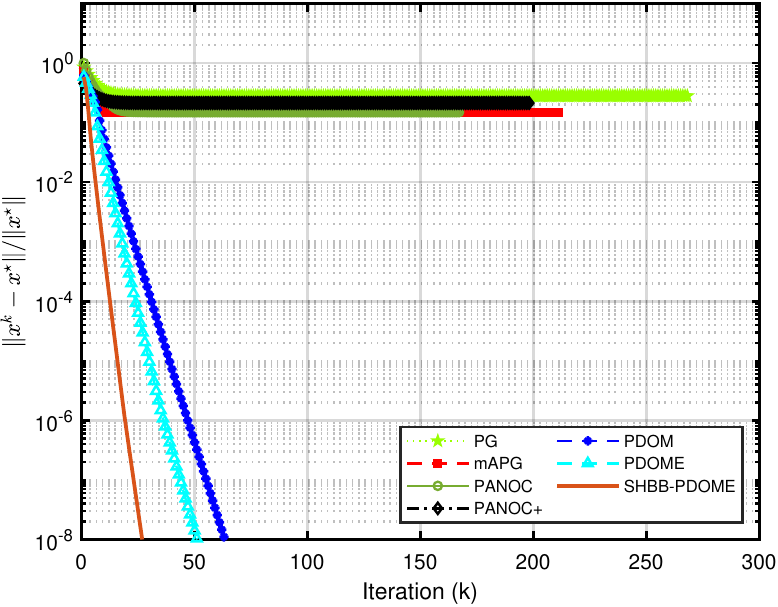}\\
    (b)
  \end{minipage}\hfill
  \begin{minipage}[t]{0.32\textwidth}
    \centering
    \includegraphics[width=\linewidth]{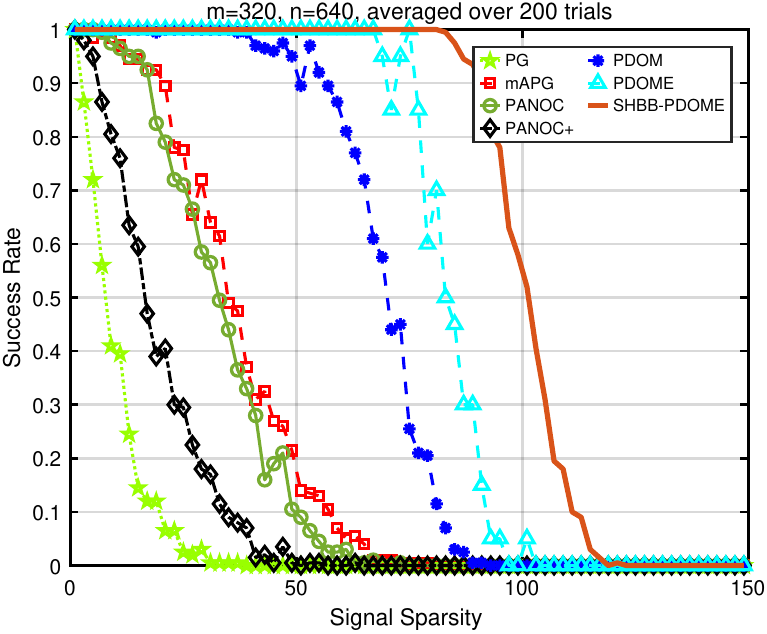}\\
    (c)
  \end{minipage}
  \caption{$\ell_0$ sparse recovery. Left: 
  performance comparison of the data residual versus the iteration count $k$.
  Middle:Normalized recovery error $\|x^k-x^\star\|/\|x^\star\|$ versus iteration count $k$.   Right: phase transition curve of $\ell_0$ sparse recovery at varying sparsities.}
  \label{NEfig2}
\end{figure}

\begin{table}[htbp]
  \centering
  \caption{Average final NRE and iteration count over 20
  independent trials for each value of $m$.
  The stopping criterion is
  $\|x^{k+1}-x^k\|/(1+\|x^{k+1}\|)<10^{-14}$.}
  \label{tab:II}
  \begin{tabular}{l cc cc cc}
    \toprule
    & \multicolumn{2}{c}{$m=100$}
    & \multicolumn{2}{c}{$m=500$}
    & \multicolumn{2}{c}{$m=1000$} \\
    \cmidrule(lr){2-3}
    \cmidrule(lr){4-5}
    \cmidrule(lr){6-7}
    Algorithm
    & NRE & Iterations
    & NRE & Iterations
    & NRE & Iterations \\
    \midrule
    PG
    & 1.06e-13 & 185.2
    & 2.57e-02 & 195.5
    & 7.50e-02 & 201.3 \\
    mAPG
    & 5.02e-13 & 119.6
    & 1.04e-13 & 141.8
    & 7.43e-14 & 148.4 \\
    PANOC
    & 7.23e-14 & 107.0
    & 5.50e-14 & 113.8
    & 5.20e-14 & 117.8 \\
    PANOC+
    & 7.28e-14 & 133.5
    & 5.62e-14 & 141.9
    & 3.11e-02 & 145.9 \\
    PDOM
    & 2.11e-14 & 58.2
    & 1.75e-14 & 62.5
    & 1.75e-14 & 64.8 \\
    PDOME
    & 1.69e-14 & 48.9
    & 1.49e-14 & 52.6
    & 1.61e-14 & 54.4 \\
    \textbf{SHBB-PDOME}
    & \textbf{3.40e-15} & \textbf{21.4}
    & \textbf{3.52e-15} & \textbf{24.2}
    & \textbf{4.07e-15} & \textbf{25.8} \\
    \bottomrule
  \end{tabular}
\end{table}
The first two subfigures (a) and (b) of Figure~\ref{NEfig2} illustrate the convergence behavior of the $\ell_0$ instances with $m=500$ from Table~\ref{tab:II}, averaged over three random seeds with shaded $\pm$ standard deviation. Compared with the benchmark algorithms, SHBB-PDOME reaches the stopping criterion in fewer iterations and achieves a smaller final NRE in these instances. Subfigure (c) of Figure~\ref{NEfig2} shows that the successful recovery rates of SHBB-PDOME  are significantly higher than those of the benchmark algorithms.

Table~\ref{tab:II} reports the average performance of SHBB-PDOME and the benchmark methods across different problem sizes in the noiseless
setting. SHBB-PDOME requires fewer iterations to satisfy
the stopping criterion and achieves a smaller final NRE
in these instances.
These results demonstrate its effectiveness
in terms of iteration counts and recovery accuracy.

\subsection{Nonconvex Sparse Approximation Problem}
Following the $\ell_{1/2}$ regularization framework in
\cite{xu2012l_}, we consider the nonconvex sparse
least-squares problem
\begin{equation}\label{eq:45}
\min_{x\in\mathbb R^n}\;
\frac12\|Ax-y\|^2+\lambda_2\|x\|_{1/2}^{1/2}.
\end{equation}
Here $A\in\mathbb R^{m\times n}$ is a Gaussian sensing matrix
with independent entries distributed as $\mathcal N(0,1/m)$,
and $y\in\mathbb R^m$ is the observation vector.
For each instance, we compute $L_s=\|A\|_2^2$ and set
$\lambda_2=0.05\|A^\top y\|_\infty$.
The smooth term $s(x)=\frac12\|Ax-y\|^2$ has the constant
positive semidefinite Hessian $\mathbf M=A^\top A$.
For $\beta>0$, the matrix
$\mathbf B_\beta=A^\top A+\beta I$ is positive definite.
We compute its Cholesky factorization once per instance
and reuse the factors to solve $\mathbf B_\beta d_N=-g_k$
at each iteration. This direct linear solve implements
the regularized Newton endpoint defined in Section~2.2,
without explicitly forming the inverse or modifying
the objective.

The $\ell_{1/2}$ quasi-norm is
$\|x\|_{1/2}=(\sum_{i=1}^{n}|x_i|^{1/2})^2$, so the sparsity
penalty $\|x\|_{1/2}^{1/2}=\sum_{i=1}^{n}|x_i|^{1/2}$
is separable. For the fixed positive definite diagonal metric
$\mathbf M_{\mathrm d}=\operatorname{Diag}(m_1,\ldots,m_n)$,
the proximal mapping can therefore be evaluated coordinatewise
using the half-thresholding formula
\cite[Theorems~1--2]{xu2012l_}.
For a steplength $a>0$ and a proximal center $c$,
the $i$th coordinate is obtained by minimizing
$\lambda_2|u|^{1/2}+\frac{m_i}{2a}(u-c_i)^2$
over $u\in\mathbb R$.
Setting $\nu_i:=a\lambda_2/m_i$, we obtain
\[
x_i^{\mathrm{prox}}=
\begin{cases}
0,
& |c_i|\le \frac32\nu_i^{2/3},\\[4pt]
\displaystyle
\frac{2c_i}{3}
\left[
1+\cos\left(
\frac{2\pi}{3}
-\frac23\arccos\left(
\frac{3\sqrt{3}\,\nu_i}{4|c_i|^{3/2}}
\right)
\right)
\right],
& |c_i|>\frac32\nu_i^{2/3}.
\end{cases}
\]
Here $a=\gamma\eta_\mu$ for the dogleg candidate and $a=\eta$
for the safeguard candidate. At the threshold, we select zero.
When $m_i=1$, the formula for coordinate $i$ reduces to the
Euclidean half-thresholding formula.

For SHBB-PDOME, we set
$\gamma=0.995$, $\zeta=0.045$, $w=0.7$,
$\rho_0=0.5$, $\delta_\rho=0.05$,
$\alpha_{\mathrm d}=0.15$,
$\beta_{\mathrm d}=10^{-2}$,
$\beta=10^{-6}$,
and
$\eta=0.99/L_s$,
which satisfies
$\eta<\lambda_{\min}^{\hat h}/L_s$
whenever $\lambda_{\min}^{\hat h}>0.99$;
this is an instance condition, verified for
the instances used below.

\begin{figure}[htbp]
    \centering
    \begin{minipage}[t]{0.32\textwidth}
        \centering
        \includegraphics[width=\linewidth]{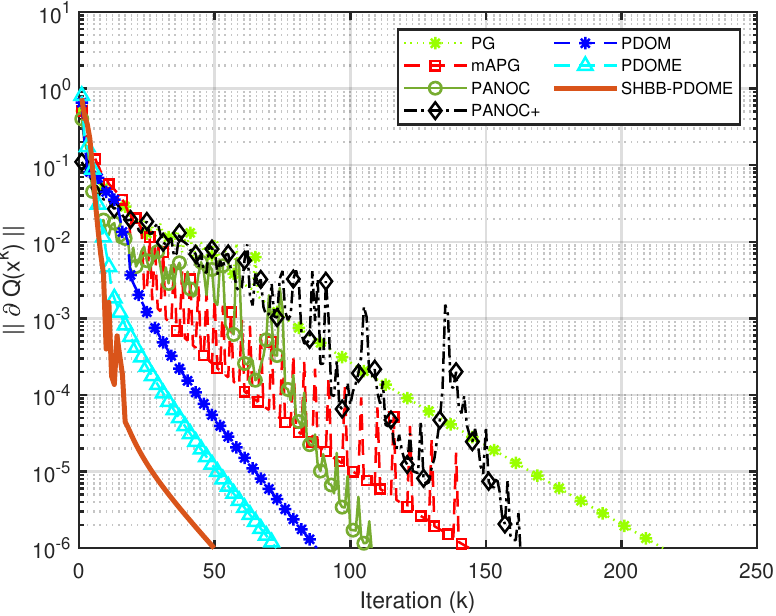}\\
    (d)
    \end{minipage}\hfill
    \begin{minipage}[t]{0.32\textwidth}
        \centering
        \includegraphics[width=\linewidth]{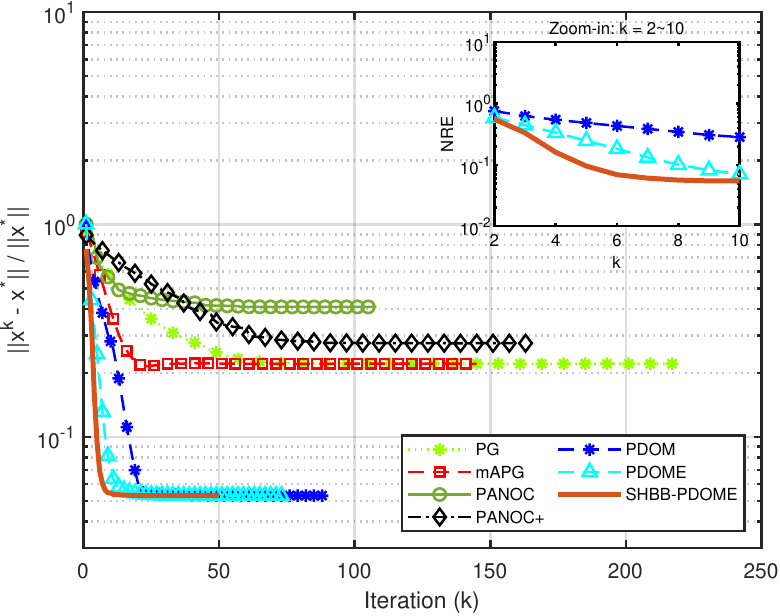}\\
    (e)
    \end{minipage}\hfill
    \begin{minipage}[t]{0.32\textwidth}
        \centering
        \includegraphics[width=\linewidth]{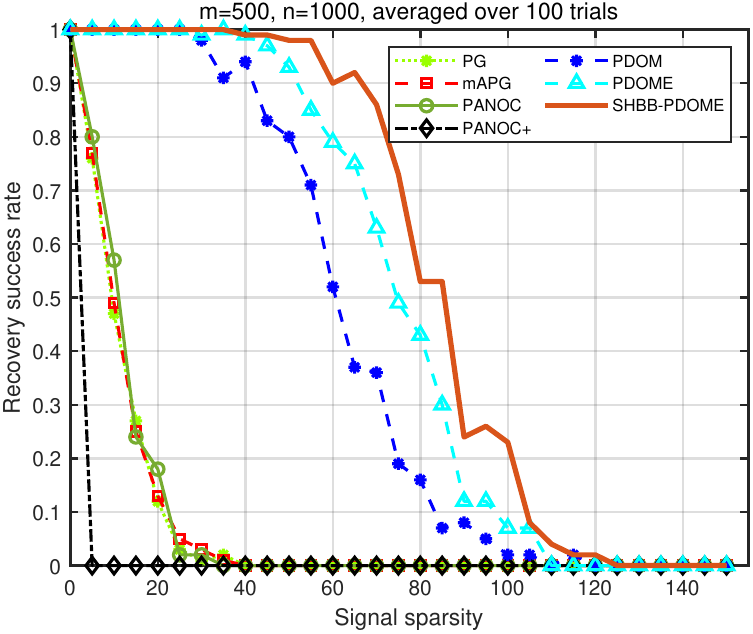}\\
    (f)
    \end{minipage}
    \caption{$\ell_{1/2}$ sparse approximation: Left: residual versus iteration count $k$, for the $\ell_{1/2}$ instances. Middle:NRE $\|x^k-x^\star\|/\|x^\star\|$ versus iteration count $k$; Right: phase transition: success rate versus the number of nonzero entries of $x^\star$.
}
    \label{NEfig3}
\end{figure}

\begin{table}[H]
  \centering
  \caption{Average final NRE and iteration count over 50
  independent trials for each value of $m$.
  The stopping criterion is
  $\|x^{k+1}-x^k\|/(1+\|x^{k+1}\|)<10^{-6}$.}
  \label{tab:III}
  \begin{tabular}{l cc cc cc}
    \toprule
    & \multicolumn{2}{c}{$m=100$}
    & \multicolumn{2}{c}{$m=500$}
    & \multicolumn{2}{c}{$m=1000$} \\
    \cmidrule(lr){2-3}
    \cmidrule(lr){4-5}
    \cmidrule(lr){6-7}
    Algorithm
    & NRE & Iterations
    & NRE & Iterations
    & NRE & Iterations \\
    \midrule
    PG
    & 7.62e-01 & 2001.0
    & 1.18e-02 & 1765.5
    & 5.14e-04 & 353.9 \\
    mAPG
    & 8.57e-02 & 741.0
    & 4.79e-04 & 207.8
    & 5.14e-04 & 110.8 \\
    PANOC
    & 9.05e-02 & 476.6
    & 4.78e-04 & 130.6
    & 5.14e-04 & 58.8 \\
    PANOC+
    & 3.61e-01 & 176.5
    & 8.96e-02 & 72.9
    & 1.05e-01 & 70.3 \\
    PDOM
    & 3.95e-04 & 159.4
    & 4.79e-04 & 56.3
    & 5.14e-04 & 44.5 \\
    PDOME
    & 5.33e-04 & 45.5
    & 5.36e-04 & 29.1
    & 5.44e-04 & 26.8 \\
    \textbf{SHBB-PDOME}
    & \textbf{5.02e-04} & \textbf{31.8}
    & \textbf{5.20e-04} & \textbf{16.6}
    & \textbf{5.34e-04} & \textbf{16.6} \\
    \bottomrule
  \end{tabular}
\end{table}

The first two subfigures (d) and (e) of Figure~\ref{NEfig3} illustrate the convergence behavior of the $\ell_{1/2}$ instances with $m=500$ from Table~\ref{tab:III}, averaged over three random seeds with shaded $\pm$ standard deviation. Compared with the benchmark algorithms, SHBB-PDOME reaches the stopping criterion in fewer iterations, while its final NRE is comparable to that of the best-performing methods. Subfigure (f) of Figure~\ref{NEfig3} shows that the successful recovery rates of SHBB-PDOME  are significantly higher than those of the benchmark algorithms.

Table~\ref{tab:III} reports the average performance of all tested algorithms on the nonconvex sparse
recovery problem across different problem sizes. The results show that SHBB-PDOME requires fewer iterations to satisfy the stopping criterion
in most settings, while its final NRE is comparable to that
of the competing methods.

\subsection{Robust Principal Component Analysis(RPCA)}
\label{sec:rpca}

We consider RPCA, which
decomposes an observed matrix $Y\in\mathbb R^{m\times n}$ into
a low-rank component $L$ and a sparse component $S$. Following
the experimental framework of \cite{candes2011robust}, we solve
\begin{equation}\label{eq:rpca}
\min_{L,S\in\mathbb R^{m\times n}}
\frac12\|Y-L-S\|_F^2
+\delta_{\mathcal R_r}(L)+\lambda_3\|S\|_0,
\qquad
\mathcal R_r:=\{L:\operatorname{rank}(L)\le r\},
\end{equation}
where $\delta_{\mathcal R_r}$ is the indicator function of
$\mathcal R_r$, $\lambda_3>0$, and $\|S\|_0$ counts the nonzero entries of $S$.

For the smooth term $s(L,S)=\frac12\|Y-L-S\|_F^2$,
the gradient Lipschitz constant is $L_s=2$ under the product
Frobenius norm. The regularized Newton direction is computed
as described in Section~\ref{sec:hybrid}.
Within the diagonal-metric framework, we set
$\mathbf M_{\mathrm d}=I$ for this experiment.
The proximal subproblem therefore uses the Frobenius metric
and, by separability of the penalty in $L$ and $S$, reduces
to a best rank-$r$ approximation and entrywise hard thresholding.
For a proximal center $(C_L,C_S)$ and a steplength $a>0$,
the updates are
\begin{equation}\label{eq:rpca-prox}
\begin{aligned}
L^+ &\in
\operatorname*{argmin}_{\operatorname{rank}(L)\le r}
\|L-C_L\|_F^2,\\
S^+_{ij}
&=
\begin{cases}
(C_S)_{ij}, & |(C_S)_{ij}|>\sqrt{2a\lambda},\\
0, & |(C_S)_{ij}|\le\sqrt{2a\lambda}.
\end{cases}
\end{aligned}
\end{equation}
The low-rank projection retains the $r$ largest singular values
of $C_L$. At the hard-thresholding boundary, we select zero.
Here $a=\gamma\eta_\mu$ for a dogleg candidate and $a=\eta$
for the safeguard candidate.

For the experiments in Table~\ref{tab:rpca_results}, we set $m=n$,
$r=\max\{1,\operatorname{round}(0.05m)\}$, and
$\lambda_3=1/\sqrt n$. The ground-truth low-rank matrix is
generated as $L^\star=UV^\top$, where
$U\in\mathbb R^{m\times r}$ and $V\in\mathbb R^{n\times r}$
have independent standard normal entries, following the
generation procedure in \cite{candes2011robust}. The sparse component
$S^\star$ contains $0.1m^2$ nonzero entries at uniformly sampled
distinct locations. Each nonzero entry equals $-1$ or $+1$
with equal probability. The observations are noiseless:
$Y=L^\star+S^\star$. We initialize $L^0$ as a best rank-$r$
approximation of $Y$, computed by truncated SVD, and set $S^0=0$,
ensuring $\operatorname{rank}(L^0)\le r$.
The results in Table~\ref{tab:rpca_results} are averaged over
20 independent trials.
\begin{figure}[htbp]
  \centering
  \begin{minipage}[t]{0.32\textwidth}
    \centering
    \includegraphics[width=\linewidth]{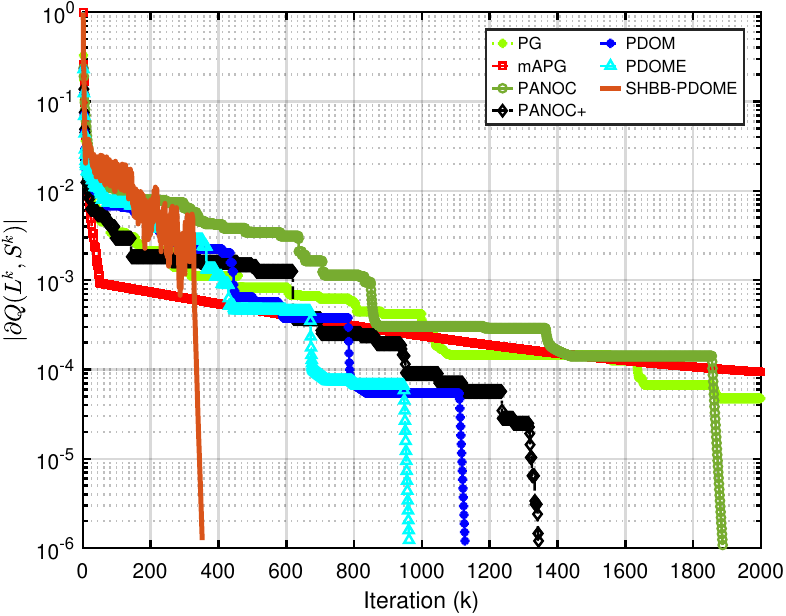}\\
    (g)
  \end{minipage}\hfill
  \begin{minipage}[t]{0.32\textwidth}
    \centering
    \includegraphics[width=\linewidth]{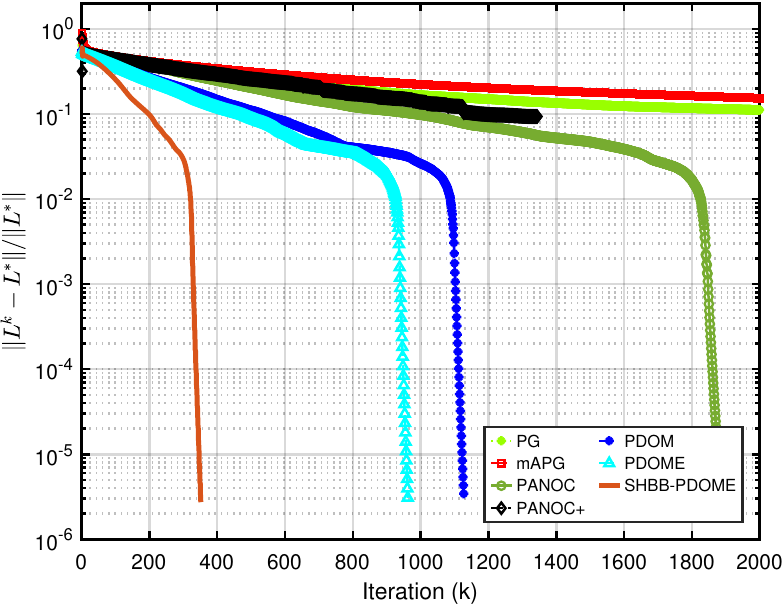}\\
    (h)
  \end{minipage}\hfill
  \begin{minipage}[t]{0.32\textwidth}
    \centering
    \includegraphics[width=\linewidth]{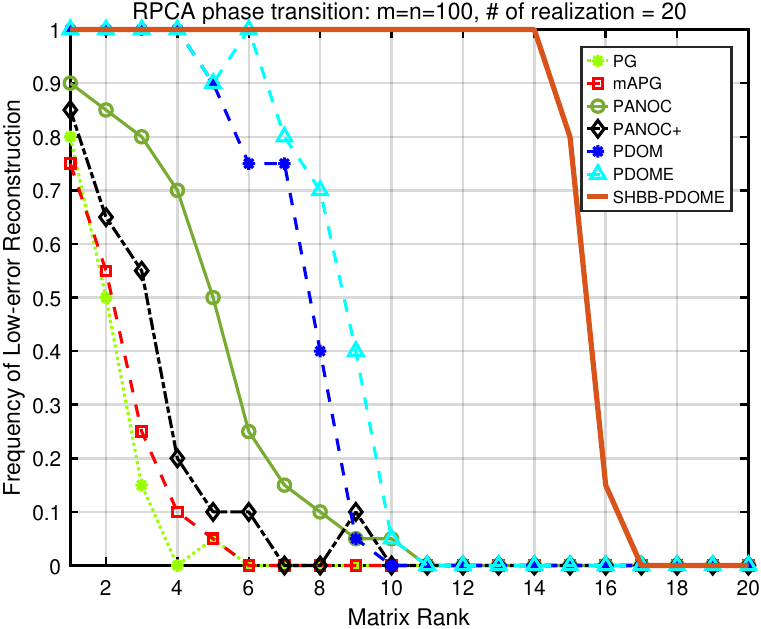}\\
    (i)
  \end{minipage}
  \caption{RPCA. Left: stationarity residual
$\operatorname{dist}(0,\partial Q(L^k,S^k))$
versus iteration count $k$.
  Middle: low-rank recovery error
  $\|L^k-L^\star\|_F/\|L^\star\|_F$ versus iteration count $k$.
  Right: frequency of low-error reconstruction versus matrix
  rank for $m=n=100$, averaged over 20 realizations.}
  \label{fig:rpca}
\end{figure}

\begin{table}[htbp]
  \centering
  \caption{Average final low-rank recovery error and iteration count
over 10 independent trials for each matrix size. The stopping criterion is
$\|x^{k+1}-x^k\|/(1+\|x^{k+1}\|)<10^{-6}$,
where $x^k=(\operatorname{vec}(L^k),\operatorname{vec}(S^k))$.}
  \label{tab:rpca_results}
  \begin{tabular}{l cc cc cc}
    \toprule
    & \multicolumn{2}{c}{$m=100$}
    & \multicolumn{2}{c}{$m=500$}
    & \multicolumn{2}{c}{$m=1000$} \\
    \cmidrule(lr){2-3}
    \cmidrule(lr){4-5}
    \cmidrule(lr){6-7}
    Algorithm
    & NRE & Iterations
    & NRE & Iterations
    & NRE & Iterations \\
    \midrule
    PG
    & 1.2996e-01 & $>2000$
    & 2.2276e-01 & $>2000$
    & 2.6002e-01 & $>2000$ \\
    mAPG
    & 1.1469e-01 & 1973
    & 2.5146e-01 & $>2000$
    & 2.9234e-01 & $>2000$ \\
    PANOC
    & 3.5619e-02 & 1982
    & 1.0344e-01 & $>2000$
    & 1.2261e-01 & $>2000$ \\
    PANOC+
    & 1.0210e-01 & 1419
    & 2.4724e-01 & $>2000$
    & 3.1294e-01 & $>2000$ \\
    PDOM
    & 4.0751e-02 & 1440
    & 6.0301e-02 & $>2000$
    & 8.3429e-02 & $>2000$ \\
    PDOME
    & 2.5298e-02 & 1225
    & 4.8262e-02 & $>2000$
    & 6.8941e-02 & $>2000$ \\
    \textbf{SHBB-PDOME}
    & \textbf{1.5461e-05} & \textbf{551}
    & \textbf{6.0747e-03} & \textbf{1162}
    & \textbf{7.7933e-06} & \textbf{817} \\
    \bottomrule
  \end{tabular}
\end{table}
For SHBB-PDOME, we use $\beta=10^{-6}$, $\gamma=0.998$,
$\zeta=10^{-6}$, $\eta=0.99/L_s$, and
$\tau_{\min}=0.01/L_s$, where $L_s=2$.
These values satisfy the extrapolation condition in
Lemma~\ref{lem:lyap}.
The BB mixing weight is $w=0.6$.
The tolerance for the inner stopping criterion is $10^{-8}$.

The first two panels of Figure~\ref{fig:rpca} show the convergence curves for one instance with $m=100$. The final panel reports recovery success rates as the rank
varies, using $0.2m^2$ nonzero entries in $S^\star$. A trial is successful if
\[
\frac{\|\widehat L-L^\star\|_F}{\|L^\star\|_F}<8\times 10^{-2},
\]
where $\widehat L$ is the recovered low-rank matrix. Table~\ref{tab:rpca_results} reports the average performance
of all tested algorithms on the RPCA problem across different
matrix sizes.

\section{Conclusion}\label{sec:5}
This paper presents SHBB-PDOME for minimizing the sum of a smooth quadratic function and a proximable, nonconvex, nonsmooth term. Within the PDOME framework, the method
constructs the dogleg path using an adaptive steplength obtained by blending a structured Barzilai--Borwein steplength
with the local steplength stored at the preceding iteration and then clipping the result. A fixed diagonal metric is used
in both the weighted BB formulas and the Bregman proximal subproblems. Under the stated assumptions, we establish
sufficient decrease of a Lyapunov function and show that every cluster point in the interior of the kernel domain
is critical. For the diagonal quadratic kernels used in the experiments, this conclusion applies to every cluster point.
Under additional assumptions, including the K{\L} property, we establish finite length of the generated sequence and
its convergence to a critical point. With an additional relative error condition, we obtain convergence rates determined by the K{\L} exponent. Numerical experiments on
nonconvex sparse recovery, sparse approximation, and robust principal component analysis demonstrate the effectiveness
of SHBB-PDOME, with faster convergence and improved recovery performance in several tested settings. Future work will
investigate extensions to nonquadratic smooth terms and alternative choices of the metric.

\section*{Declarations}
\noindent{\bf Funding:} This research is supported by the National Natural Science Foundation of China (NSFC) grants 92473208, 12401415, the Key Program of National Natural Science of China 12331011, the 111 Project (No. D23017),  the Natural Science Foundation of Hunan Province (No. 2025JJ60009).

\noindent{\bf Data Availability:} Enquiries about data/code availability should be directed to the authors.

\noindent{\bf Competing interests:} The authors have no competing interests to declare that are relevant to the content of this paper.
\bibliographystyle{elsarticle-num}	
\bibliography{cas-refs.bib}

% %% For authoryear reference style

% \bibitem[Lamport(1994)]{lamport94}
%   \textit{\LaTeX: a document preparation system},
%   Addison Wesley, Massachusetts,

% \end{thebibliography}
\end{document}